\documentclass[journal,draftcls,onecolumn,letterpaper,twoside,11pt]{IEEEtran}

\usepackage{amsmath}
\usepackage{subfig}
\usepackage{algorithm,algorithmic}
\usepackage{dsfont}
\usepackage{amssymb}
\newtheorem{theorem}{Theorem}
\newtheorem{corollary}{Corollary}[theorem]
\newtheorem{lemma}[theorem]{Lemma}

\newcommand{\bR}{\mathbb{R}}
\newcommand{\bN}{\mathbb{N}}
\newcommand{\bs}{\boldsymbol}
\newcommand{\bS}{\mathbb{S}}
\newcommand{\cS}{\mathcal{S}}
\newcommand{\cJ}{\mathcal{J}}

\newcommand{\cC}{\mathcal{C}}
\newcommand{\cN}{\mathcal{N}}
\newcommand{\cL}{\mathcal{L}}
\newcommand{\cY}{\mathcal{Y}}
\newcommand{\cF}{\mathcal{F}}
\DeclareMathOperator*{\esssup}{esssup}
\newcommand{\Pro}{{\mathbb{P}}}
\newcommand{\Exp}{{\mathbb{E}}}
\newcommand{\sD}{{\mathsf{Div}}}
\DeclareMathOperator*{\argmax}{arg\,max}
\DeclareMathOperator*{\argmin}{arg\,min}
\usepackage[authoryear]{natbib}
\usepackage{graphicx}
\graphicspath{{./Images/}}

\begin{document}

\title{Sequential Change Diagnosis Revisited \\
and the Adaptive Matrix CuSum}
%


\IEEEauthorblockN{
\author{\IEEEauthorblockN{Austin Warner and Georgios Fellouris}\\
\IEEEauthorblockA{\textit{Department of Statistics} \\
\textit{University of Illinois at Urbana-Champaign}\\
Urbana, USA \\
$\lbrace$ awarner5, fellouri $\rbrace$@illinois.edu}
\thanks{This research was supported by the US National Science Foundation under
grant AMPS 1736454 through the University of Illinois at Urbana-Champaign.}}}


        \maketitle



\begin{abstract}
The problem of sequential change diagnosis is considered, where observations are obtained  on-line, an  abrupt change  occurs in their distribution, and the goal is to  quickly detect the change  and accurately identify  the post-change distribution, while controlling the false alarm rate. A finite set of alternatives is postulated for the post-change regime, but no prior information  is assumed for  the unknown change-point.  A  drawback of many algorithms that have been proposed for this problem  is the  implicit use of  pre-change data for determining the post-change distribution. This  can lead to very large  conditional probabilities of misidentification, given that there was no false alarm,  unless the change occurs soon after monitoring begins. A novel, recursive  algorithm is proposed and shown to resolve this issue without the use of additional tuning parameters and without sacrificing control of the worst-case  delay  in Lorden's  sense. A theoretical analysis is conducted for a general family of sequential change diagnosis procedures, which  supports the proposed algorithm and  revises certain state-of-the-art results.  Additionally, a  novel, comprehensive method is proposed for the design and evaluation of sequential change diagnosis algorithms. This method is illustrated with  simulation studies, where  existing procedures are compared to the proposed.
\end{abstract}




\section{Introduction}
The problem of quickly  detecting a change in  sequentially acquired data dates back to the  pioneering works \cite{shewhart1931economic} and \cite{page1954continuous}  and is motivated by a wide range of engineering and scientific applications. Examples of 
 such applications can be found in  industrial process quality control~(\cite{bissell1969cusum,hawkins2003changepoint,joe2003statistical}), target detection and identification (\cite{ru2009detection,blackman2004multiple}),  integrity monitoring of navigation systems (\cite{nikiforov1993application,bakhache1999reliable}),
 target tracking (\cite{tartakovsky2003sequential}), network intrusion detection (\cite{tartakovsky2006detection,tartakovsky2006novel}), bioterrorism (\cite{rolka2007issues,fienberg2005statistical}), and genomics (\cite{siegmund2013change}).  In most of these applications,  there are  many  possible 
types of change  and  it is useful, if not critical,  to not only detect the change quickly but also to correctly identify  its  type upon stopping.  

The problem of simultaneously detecting a change in the distribution of  sequentially collected data and  identifying  the correct post-change distribution among a finite set of   alternatives is  known as \textit{sequential change diagnosis}. The literature on this problem  can be broadly classified into two categories.  In the first one (see, e.g., \cite{dayanik2008bayesian,dayanik2013asymptotically, ma2021two,tartakovsky2021asymptotic}),    a Bayesian formulation is adopted that assumes a prior distribution on the change-point,  and sometimes on the type of change. In the  second (see, e.g.,  \cite{nikiforov1995generalized, tartakovsky2008multidecision, lai2000sequential, nikiforov2000simple, pergamenchtchikov2022minimax, oskiper2002online}),   no such prior information is assumed. Both setups are considered in some works, such as \cite{lai2000sequential, pergamenchtchikov2022minimax}.

A Bayesian formulation of the  sequential change \textit{detection} problem was first proposed in  \cite{shiryaev1963optimum} (see also \cite{shiryaev2007optimal}), where the Bayes rule was derived under the assumptions of a Geometric prior for the change-point, iid observations  before and after the change conditionally on the change-point, and completely specified distributions.   A Bayesian formulation of the sequential change \textit{diagnosis} problem  was  first proposed in \cite{malladi1999generalized} and \cite{dayanik2008bayesian}. In the latter work,  given a Geometric prior distribution for the change-point, the Bayes rule is shown to raise an alarm the first time the multi-dimensional posterior probability process enters  a union of disjoint regions, each of which corresponds to a posited post-change distribution.  These  regions need to be computed off-line, via dynamic programming, and the necessary calculations can be very demanding.
 To address this issue, a computationally feasible  procedure was  proposed in \cite{dayanik2013asymptotically} and was shown to achieve the Bayes risk asymptotically
 as the probabilities of false alarm and false identification go to 0.  A two-stage approach  was considered in \cite{ma2021two}, where a sequential change detection algorithm is applied and, after it is declared that the change has occurred, a  sequential hypothesis test is then performed  to determine the post-change distribution. Non-iid models for the pre-change and post-change regimes, as well as a non-geometric prior for the change-point, were  allowed in \cite{tartakovsky2021asymptotic} in a special case of the \textit{multichannel} problem. The latter refers to a case of the sequential diagnosis problem wherein a number of data channels are observed in parallel and a change occurs in the distribution of an unknown  subset of the channels. In \cite{tartakovsky2021asymptotic}, the channels are assumed to be independent of each other and it is a priori known that the change occurs in  only one of them.

In the absence of a priori information regarding the change-point, the sequential change diagnosis problem turns out to be significantly more complex than the pure sequential change detection  problem, even in the case of  iid observations before and after the change. Indeed, a recursive algorithm for the sequential change detection problem, Page's Cumulative Sum (CuSum), was proposed early on in \cite{page1954continuous} and, since then,  has been  grounded on strong theoretical support. In particular, assuming completely specified pre-change and post-change distributions,  it was shown in  \cite{moustakides1986optimal}  to  minimize  \textit{Lorden's  criterion} (\cite{lorden1971procedures}), i.e.,  the worst-case conditional expected detection delay with respect to both the change-point and the history of the observations up 
to the change-point,  subject to a user-specified bound on the false alarm rate.   Modifications of  the CuSum algorithm  have  been introduced to allow for parametric uncertainty and/or temporal dependence in the pre-change and post-change models, e.g., \cite{mei2006sequential, siegmund1995using, lai1998information, lai2010sequential, lai1999efficient}.

 The first sequential change diagnosis scheme that does not  utilize prior information regarding the change-point (and in general) was
 proposed  in \cite{nikiforov1995generalized}, and was a  generalization of the CuSum algorithm, termed the ``Generalized CuSum.'' However, unlike the CuSum,  this algorithm does not  admit a recursive computational structure,  even in the iid setting, and the number of operations it requires per time instance grows with the number of collected observations. Since the change can take a very long time to occur,  this procedure needs to be modified to be applicable in practice. One such modification, suggested in   \cite{lai2000sequential}, uses at any given time only the most recent data via a sliding window of deterministic size, and is referred to as the ``window-limited Generalized CuSum.'' A different modification   of the Generalized CuSum was  proposed in  \cite{oskiper2002online} (see also \cite{tartakovsky2008multidecision}), where the order of maximizing over possible change-points and minimizing over post-change scenarios is reversed. This procedure was termed the ``Matrix CuSum,'' as it  requires the  parallel computation of a matrix of CuSum statistics.

Independently of their very different  computational requirements,   the Generalized CuSum, its window-limited modification,  and the Matrix CuSum all  control Lorden's  worst-case criterion for the  detection delay and have been shown theoretically to control  false identification metrics  \textit{only when the change occurs from the beginning of monitoring}. This, however, may be  the \textit{best} possible scenario  for these algorithms with respect to the identification task.  In fact, the  Generalized CuSum  and the Matrix CuSum have both been reported   to have, in certain cases, large probabilities of false isolation when the change-point is large \cite[Section 4.3.1]{nikiforov2016sequential}. Despite this, mainly due to its practically convenient,  recursive structure, the Matrix CuSum continues to inspire formulations of the sequential change diagnosis problem (see, e.g., \cite{huang2021asymptotic})   in which   only this
``easy'' case is considered for  the false isolation metric. 

 Intuitively, the identification ability of the Generalized CuSum and the Matrix CuSum   is compromised 
by the impact of data from before the change in estimating the post-change distribution. This impact is reduced when applying the window-limited Generalized CuSum  in  \cite{lai2000sequential}. However, the  size of the deterministic window is a tuning parameter, which is  not determined by the user-specified error constraints and whose choice   is  characterized by a fundamental trade-off. Indeed, the window should be large enough for the algorithm to be able to detect the change, especially when the actual post-change alternative is ``close'' to the pre-change distribution,  but not too large to compromise its identification ability.  An asymptotic lower bound is obtained for this window  in  \cite{lai2000sequential}, but to the best of our knowledge it is not clear how to precisely select this window, in general, in order to simultaneously guarantee the desired isolation control.

A recursive algorithm  that does not require any tuning parameters was proposed and  analyzed in  \cite{nikiforov2000simple}. This algorithm  requires the computation of the  CuSum statistics that compare the post-change distributions against the pre-change, and was termed the ``Vector CuSum,'' as its implementation requires the recursive calculation of a vector of CuSum statistics. It was argued that it can control the conditional probability of  a wrong identification, given that there was no false alarm, uniformly with respect to the change-point. However, its detection delay analysis was conducted   using  Pollak's  criterion   (\cite{pollak1985optimal}), not Lorden's.  That is, the Vector CuSum was shown to control  the worst-case conditional expected detection delay with respect to the change-point, but not also the data before the change.    In  \cite{lai2000sequential} it was  shown that a window-limited algorithm modification of the Vector CuSum
can control the probability of false identification within a fixed window from the change-point, uniformly with respect to the change-point.  However, in addition to requiring the specification of a tuning parameter,  this algorithm was not shown to control  a worst-case criterion for the detection delay, but only weighted averages of the expected detection delay.

Finally, more recently, an algorithm was proposed in  \cite{pergamenchtchikov2022minimax} that controls the conditional probability of a false identification, given that there was no false alarm,  in the special case of the multichannel problem where the channels are independent and  the change can occur in only one  of them. However, the conditional probability of false identification is only controlled when the change occurs within a fixed window of the outset. Similarly to the above window-limited algorithms,   this procedure requires the specification of certain tuning parameters,  in addition to the thresholds that guarantee the prescribed error constraints.   Similarly to the Vector CuSum, it was not shown to control's  Lorden's criterion for the delay, but a Pollak-type delay criterion. \\


The first   goal  of this paper is to  introduce  and analyze a novel sequential change diagnosis algorithm  that (i) can provide the desired guarantees for false isolation control, (ii)   does not require the specification of tuning parameters, (iii) has a practically convenient, recursive structure, and (iv) controls Lorden's criterion for the  detection delay. Specifically, the proposed scheme  is a modification of the Matrix CuSum, in which each of its statistics is reset at some  \textit{random, data-dependent} time that corresponds to an estimate of the change-point and, for this reason, we refer to it as the ``Adaptive Matrix CuSum.''
These adaptive resets  suppress  the inclusion of pre-change data when determining the post-change regime, and, in this way, they  lead to the desired false identification error control  without sacrificing the detection ability of the algorithm.  
%

The second goal  of the present work is to establish a  general theory for a family of sequential change diagnosis procedures, which  encompasses  the Vector CuSum,  the Matrix CuSum, and
its proposed modification. In particular, we establish  a uniform exponential bound for the  tail probabilities under the pre-change distribution  of the statistics of  the Vector CuSum and of the proposed scheme. This  implies a novel  lower bound on their expected time to false alarm, and suggests that pre-change data do not influence the decisions of these two schemes  about the post-change distribution.  However, we argue that, based on the techniques used in the present paper as well as in \cite{nikiforov2000simple},  the false identification control of the Vector CuSum  can  indeed be established for a certain  family of change-points, but not necessarily for all possible change-points, as reported in \cite{nikiforov2000simple}.

A final goal of the present paper is to develop a comprehensive, computationally efficient  framework for  the design and fair comparison of sequential change diagnosis schemes. Based on this framework, we conduct a number of simulation studies in which we compare various sequential change diagnosis schemes under two different setups of the multichannel problem.

The remainder of the paper is organized as follows: In Section \ref{sec:state} we formulate mathematically the problem of sequential change diagnosis  that we consider in this work. 
In Section~\ref{sec:CUSUM} 
we review the CuSum statistic and the min-CuSum procedure.
In Section~\ref{sec:family} we introduce a family of sequential change diagnosis procedures  and  the proposed novel procedure. In Section~\ref{sec:suppressing} we provide a description of how, for some algorithms, pre-change data can systematically lead to an  erroneous decision regarding the post-change distribution.  In 
Section~\ref{sec:analysis} we analyze the performance of  schemes in the family which we introduce,  with an emphasis on the proposed procedure. In Section~\ref{sec:design} we introduce the proposed framework for the design and comparison of sequential change diagnosis  procedures in this family, which is  followed by simulation studies in said framework in Section~\ref{sec:simulation}. Lastly, concluding remarks are provided in Section~\ref{sec:conclusion}.  Appendixes~A-F are dedicated to the  proofs of all main results in this work, as well as  to the statement and proof of supporting   lemmas, and Appendix G contains additional figures.

We end this section by introducing some notation that we use throughout the paper. Thus, $\bN$ is the set of positive integers, $\bN_0$ the set of non-negative integers, and $\bR$ the set of real numbers.  For $n \in \bN$, we set ${[n] \equiv \lbrace 1, 
\ldots, n \rbrace}$. For $x \in \bR$,  $x^+$  is its  positive part and $x^-$ its negative  part, i.e., $ x^+ \equiv \max\{x,0\}$ and  $ x^- \equiv \max\{-x,0\}$. We denote by $ x \wedge y $ the minimum of $ x, y $.  If $(x_n), (y_n)$ are two families of positive  numbers, where $n \in \bN$, then $x_n \sim y_n$ and $ x_n \gg y_n$ stand for $ x_n/y_n \to 1 $ and $ x_n - y_n \to \infty $, respectively, as $ n \to \infty. $  
We denote by $\cN(\mu, 1)$ the density of a Gaussian random variable with mean $\mu$ and standard deviation $1$.
The indicator function is written as $ \mathds{1}(\cdot)$.

\section{Problem Statement} \label{sec:state}
Let  $X\equiv  \{X_n, n \in \bN\}$  be a  sequence of  independent \mbox{$\bS$-valued}  random elements, where $\bS$ is an arbitrary Polish space. We assume that   each term of $X$  has  a  density  with respect to a $\sigma$-finite measure $\lambda$, which  is  $f$ up to and including  some deterministic time $\nu \in  \bN_0 $ and is $g$ after $\nu$.  The pre-change density, $f$, is completely specified, but there is  a finite number,   $K$,  of posited  alternatives  for the post-change density, i.e., ${g \in \{g_1, \ldots, g_K\} }$, and the \textit{change-point},~$\nu$, is completely unknown. The terms of the sequence $X$ are observed sequentially, and  the goal  is to  quickly detect the change and  also  identify the correct post-change distribution  upon detection. Thus, we need to specify an $\bN$-valued random time, $T$, at which we declare that the change has occurred, and a $[K]$-valued random variable, $D$, that represents the decision regarding the post-change density at the time of stopping. That is, for any $n \in \bN$ and $i \in[K]$,    the alarm is raised and   $g_i$  is  declared as the correct post-change density after having taken $n$ observations on the event  $\{T=n, D=i \}$.  We refer to  such a pair $(T,D)$ as a \textit{sequential change diagnosis procedure} if 
 $T$ is  an $\{\cF_n, n \in \bN \}$-stopping time and   $D$  is an $\mathcal{F}_T$-measurable, $[K]$-valued random variable, or in other words if  
 $$\{T=n \}, \{T=n, D=i \} \in \cF_n \quad \text{ for all } \; n \in \bN, \, i \in [K],$$ 
where  $\cF_n$ is the $\sigma$-algebra generated by the first $n$ observations, $X_1, \ldots, X_n$, and $ \mathcal{F}_0 $ is the trivial $\sigma$-algebra. We denote by $ \mathcal{C}$ the family of all sequential change diagnosis procedures.

We denote by~$ \mathbb{P}_{\infty}$ the  distribution of $X$, and by~$\mathbb{E}_{\infty}$ the corresponding expectation, when the change never occurs, i.e., when 
$X$ is a sequence of independent  random elements with common density $f$. We denote by  $ \mathbb{P}_{\nu,i} $ the distribution of $X$, and by $ \mathbb{E}_{\nu,i} $ the corresponding expectation,  when the change occurs at time $ \nu $ and the post-change density is $g_i$. For simplicity, when the change occurs from the very beginning,  we suppress the dependence on the change-point  and  set $ \mathbb{P}_i \equiv\mathbb{P}_{0,i}$ and  $\mathbb{E}_i \equiv \mathbb{E}_{0,i}$.   Moreover, without loss of generality,  we  restrict ourselves to stopping times that are not almost surely bounded under any  $\Pro_{\nu,i}$, where  $\nu \in  \bN_0$, $i \in [K]$.

To measure the ability of  a sequential change diagnosis procedure  $(T,D) \in \cC$   to avoid  false alarms we use the average number of  observations  until stopping  under $\Pro_\infty$, i.e., $\Exp_\infty[T]$. We denote by  $\cC(\alpha)$ the subfamily of diagnosis procedures whose expected time until stopping  under $\Pro_\infty$ is at least $1/ \alpha $,  i.e., 
 \begin{equation}
 \cC(\alpha) \equiv \{(T,D) \in \cC: \Exp_{\infty}[T] \geq 1/\alpha\},
 \end{equation}
where $\alpha \in (0,1)$ is a  user-specified value that represents tolerance to false alarms. 
 
To measure the  ability of $(T,D) \in \cC$  to  isolate the  post-change density $g_i$ when the change occurs at time $\nu$, we use the conditional probability of  an incorrect identification  given that there was no false alarm, i.e., 
$ \Pro_{\nu, i}(D \neq i | T > \nu)$. Thus, we denote  by $\cC(\alpha, \beta, N)$ the subfamily of change diagnosis procedures in $\cC(\alpha)$ that control the  conditional probability  of  false isolation  below $ \beta$ when the change-point belongs to a set $N \subseteq \bN_0$, i.e., 
\begin{equation}
\cC(\alpha, \beta, N) \equiv   \big\lbrace (T,D) \in \cC(\alpha) :\max_{i \in [K]} \sup_{\nu \in N}  \Pro_{\nu, i}(D \neq i \, | \, T > \nu)\leq \beta \big\rbrace,
\end{equation}
where $\beta \in (0,1)$ is a  user-specified value that represents tolerance to  false isolations.  

Finally, to measure the  ability of $(T,D) \in \cC$ to quickly detect the change when the post-change density is $g_i$  for some $i \in [K]$,  we adopt  Lorden's  criterion  (\cite{lorden1971procedures}) by employing the worst-case conditional expected detection delay  with respect to both  the change-point and  the data up to the change:
  \begin{equation}\label{eqn:lorden_delay_def}
 \cJ_i[T]\equiv \sup_{\nu \in \bN_0} \esssup \Exp_{\nu,i}[ T - \nu \, | \,  \mathcal{F}_{\nu}, T>\nu].
\end{equation}

The problem we consider in this work is to find a     change diagnosis scheme that can be designed to belong to $\cC(\alpha, \beta, N)$ for  arbitrary  $\alpha, \beta \in (0,1)$, with  $N$ being as large as possible, and achieve
\begin{equation} \label{infimum}
 \inf_{ (T,D) \in \cC(\alpha, \beta, N)}  \cJ_i[T],
 \end{equation}
 to a first-order asymptotic approximation as $ \alpha$ and $ \beta $  go to 0, simultaneously  for every possible post-change alternative, i.e., for every $i \in [K]$.

 \subsection{Assumptions} \label{sec: notation}
Our standing assumption throughout the paper is that the  Kullback-Leibler divergences,
  \begin{align} \label{KL}
  \begin{split}
 I_{i} &\equiv   \sD(g_i  \, || \,f)  \equiv  \int \log(g_i/f) \, g_i \, d \lambda,\\
  I_{ij} &\equiv  \sD(g_i \,  || \, g_j) \equiv \int \log(g_i/g_j) \, g_i \, d \lambda,
 \end{split}
\end{align}
are positive and finite for every $i, j \in [K]$ such that $i \neq j$.   For every  $n \in \bN$ and  $i,j  \in [K]$ with $j \neq i$ we set 
\begin{equation} \label{LLR}
  \ell_{i}(n) \equiv \log  \frac{ g_i (X_n)}{ f(X_n) }, \quad \quad 
  \ell_{ij}(n) \equiv \log \frac{ g_i (X_n)}{ g_j(X_n)},
\end{equation}
so  that $I_i= \Exp_i[ \ell_i(n)]$ and $I_{ij}= \Exp_i[ \ell_{ij}(n)],$ and we  denote by  $\psi_{ij}$ the cumulant generating function of~$\ell_{ij}(1)$ under $\Pro_\infty$, i.e., 
\begin{equation} 
\psi_{ij}(\theta) \equiv \log \left( \Exp_{\infty} \left[ \exp [ \theta \, \ell_{ij}(1) ] \right] \right), \quad \theta \in \bR.
\end{equation}
For the main results of this work we need to assume that 
\begin{equation}\label{assum:A1_prime}
\psi_{ij}  \; \; \text{is finite around 0 for every} \; i,j \in [K],  \, i \neq  j,
\end{equation}
but we state this assumption explicitly when we make it.


\subsection{Example: the multichannel problem}\label{subsec:multichannel}

We illustrate the above sequential change diagnosis problem in the special case  of the
 multichannel  problem, where  $d$ independent channels are simultaneously monitored and there is a change in the marginal distributions of  an unknown subset of them at some  unknown time, $\nu$. Specifically,  suppose that channel~$i$ generates a sequence of independent $\bS_i$-valued random elements, where $\bS_i$ is some  Polish space, and  let $X_{i,n}$ denote the observation from channel $i$ at time $n$.  If the change does not occur in channel~$i$, then $X_{i,n}$ has density~$p_i$ with respect to a $\sigma$-finite measure $\lambda_i$ on $\bS_i$, whereas if the change does occur in that channel at time $\nu$, then the density of $X_{i,n}$  is~$p_i$ for  $n \leq \nu$ and $q_i$  for $n > \nu$. 
 
 In this context we have
 \begin{align} \label{multichannel_X}
 X_n &= (X_{1,n}, X_{2,n}, \ldots, X_{d,n}) \in \bS  \equiv \bS_1  \otimes \cdots \otimes \bS_d,
 \end{align} 
 and the pre-change density is
 \begin{equation}
 f(x_1, \ldots , x_d) = \prod_{i=1}^d p_i(x_{i}), \quad  \quad (x_1, \ldots, x_d) \in \bS. 
 \end{equation}
 If the  change can occur in only one channel, as it is often  assumed in the literature,   the number of posited post-change distributions,  $K$, is equal to the number of channels, $d$, and  the post-change densities are 
 \begin{equation} \label{multichannel_single}
 g_i (x_1, \ldots, x_d) = q_i(x_i) \prod_{j \in [d]: \, j \neq i} p_j(x_j), \quad \quad 
(x_1, \ldots, x_d) \in \bS, \quad i \in [d].
 \end{equation}
However, it is conceptually and practically relevant to allow for the possibility of the change occurring  in multiple channels simultaneously (\cite{mei2010efficient, xie2013sequential, fellouris2016second}). As we will see in Section \ref{sec:simulation} and elsewhere, this more general setup turns out to be much more challenging for existing change diagnosis algorithms even in the case of two channels ($d=2$),  in which case the total number of post-change alternatives is  only $K=3$, specifically, 
\begin{align}\label{multichannel_simultaneous}
  g_1(x_1, x_2) \equiv q_1 \left( x_1 \right) \cdot p_2 \left( x_{2} \right),  \quad   g_2(x_1, x_2) \equiv p_1 \left( x_{1} \right) \cdot q_2 \left( x_{2} \right), \quad
g_3(x_1, x_2) \equiv  q_1(x_1) \cdot q_2(x_2).
\end{align}

\section{The min-CuSum  Algorithm}\label{sec:CUSUM}

In this section we review  the sequential change detection algorithm that provides the basis for the change diagnosis schemes that we consider in this work.   Page's CuSum algorithm (\cite{page1954continuous}) for detecting a  change from $f$ to $g_i$, for some fixed $i \in [K]$,  raises an alarm as soon as   the statistic
\begin{align} \label{cusum_statistic}
Y_{i} (n) &\equiv \max_{ 0 \leq t \leq n}  \sum_{u=t+1}^n  \ell_{i}(u),
\quad n \in \bN,
\end{align} 
exceeds a threshold $b_i>0$,  i.e., at 
\begin{align} \label{cusum}
\begin{split}
\sigma_i (b_i) &\equiv \inf \lbrace n \geq 1 : Y_i(n) \geq b_i \rbrace,
 \end{split}
\end{align}
where we adopt the convention $ \sum_{n+1}^{n} = 0. $
An important property of this stopping rule concerning its implementation in practice is that its statistic can be computed via the following recursion:
 \begin{align}\label{cusum_recursion}
Y_{i} (n)  &= \left(Y_{i}(n-1) + \ell_{i}(n) \right)^+, \quad n \in \bN,
\end{align} 
with  $Y_{i}(0) = 0$.
Of course, this algorithm is directly applicable in our setup  only when $K=1$. When $K>1$, a standard  approach for detecting the change  is to   run in parallel the $K$ CuSum statistics, $Y_1, \ldots, Y_K$, and  stop as soon as one of them exceeds its corresponding threshold, i.e.,   at $ \min_{i \in [K]} \sigma_i(b_i).$
In what follows, we set $ b_1 = \ldots = b_K = b$ and refer to the  stopping time
\begin{equation}
 \sigma(b) \equiv \min_{i \in [K]} \sigma_i(b),
\end{equation}
as  the ``min-CuSum'' stopping time.   It is known \cite[Chapter 9.2]{tartakovsky2014sequential} that, for any $b>0$,
\begin{align} \label{CUSUM_ARL}
 \Exp_{\infty}[\sigma(b)] &\geq e^{b}/K.
\end{align} 
As a result,   $\sigma(b_\alpha) \in \cC (\alpha)$ for any $\alpha \in (0,1)$, where  
\begin{equation}\label{eqn:b_alpha}
 b_\alpha \equiv |\log \alpha| + \log K.
\end{equation}

It is also well known (see, e.g., \cite{lorden1971procedures}) that   $\sigma (b_\alpha)$ minimizes $\cJ_i$ in  $\cC (\alpha)$,  for every $i \in [K]$, to  a first-order asymptotic approximation as $\alpha \to 0$, i.e,   for every $i \in [K]$, as $\alpha \to 0$, 
 \begin{align} \label{CUSUM_AO}
  \cJ_{i} \left[ \sigma(b_\alpha) \right]  \sim  
  \frac{|\log \alpha|}{I_{i}}  \sim  \inf_{(T, D)  \in \mathcal{C}(\alpha)} \mathcal{J}_i[T].
 \end{align}

The stopping time $\sigma(b)$  is associated with a natural identification rule, $ \widehat{\sigma}(b)$, which is to select a post-change alternative  that corresponds to the largest CuSum  statistic at the time of stopping, i.e., 
\begin{equation}\label{CUSUM_isolation_rule}
 \widehat{\sigma}(b)  \in    \argmax_{i \in [K]} Y_i(\sigma(b)) ,
\end{equation}
solving ties, if any,  in some arbitrary way. The min-CuSum procedure, i.e., the pair  $(\sigma(b), \widehat{\sigma}(b))$, is often used in applications as an ad-hoc solution to  the sequential change diagnosis problem (see, e.g., \cite{chen2015quickest}), whereas 
there is some theoretical justification for its use (\cite{warnercusum, han2007detection}), which we return to  in Section \ref{sec:conclusion}.

\section{Diagnosis Procedures}\label{sec:family}

 In this section we  introduce   a novel  sequential change diagnosis scheme in the context of a more  general family that encompasses many existing procedures in the literature. 

\subsection{A family of diagnosis procedures}
For each $i \in [K]$ and $n \in \bN$, let  $W_{i}(n)$ be an \mbox{$\cF_n$-measurable} statistic,  large values of which should provide evidence that the correct post-change hypothesis is  $g_i$. Given such statistics,  for every $i \in [K]$ we denote by  $\tau_i(b,h)$ the first time  the CuSum statistic $Y_i$ is equal to or exceeds a threshold $b>0$, and, at the same time,  $W_i$  is equal to or exceeds a distinct threshold $h>0$, i.e., 
\begin{align}  \label{general_stopping_time} 
 \tau_i(b,h) &\equiv \inf \left\lbrace n \in \bN : Y_i(n) \geq b \quad \& \quad W_i(n) \geq h \right\rbrace. 
 \end{align}
Then, a natural    diagnosis procedure is
  \begin{align}  \label{family} 
 \tau(b,h) &\equiv  \min_{i \in [K]} \tau_i(b,h), \quad \widehat{\tau}(b,h)  \in   \argmin_{i \in [K]} \tau_i(b,h),
 \end{align}
 where ties in the decision rule are settled in some arbitrary way. 

With an appropriate selection of  $W_1, \ldots, W_K$, 
we can recover many of the sequential change diagnosis algorithms that have been proposed in the literature. For example, when   $W_{i}(n)$  is of the form 
\begin{align} \label{window_limited}
 \max_{ M_i(n) \leq t \leq n  } \, \min \left\lbrace \sum_{u=t+1}^n  \ell_{i}(u),  \min_{ j \in [K]: j \neq i}
\,  \sum_{u=t+1}^n  \ell_{ij}(u) \right\rbrace , 
\end{align}
we recover the Generalized CuSum in \cite{nikiforov1995generalized} when ${M_i(n)=0}$, and the window-limited Generalized CuSum in \cite{lai2000sequential}  when   ${M_i(n)=n-m}$ for some fixed $m \in \bN$.  

Other  change diagnosis schemes in the literature  can be recovered by setting  $W_i(n)$  equal to 
\begin{equation}\label{family4}
 \min_{ j \in [K]: j \neq i} W_{ij}(n), 
\end{equation}
where each  $W_{ij}(n)$ is an  $\cF_n$-measurable statistic, large values of which should provide evidence that $g_i$ is a more plausible post-change alternative than $g_j$. For example, when \begin{equation}
W_{ij}(n) = Y_i(n) - Y_j(n),
\end{equation} where $Y_i$ is the CuSum statistic defined in \eqref{cusum_statistic}, we recover   the Vector CuSum  in  \cite{nikiforov2000simple},  whereas when $W_{ij}$ is selected as the CuSum statistic for detecting a change from $g_j$ to $g_i$, i.e., 
\begin{align}\label{matrix_cusum}
Y_{ij} (n) &\equiv \max_{0 \leq t \leq n }\sum_{u=t+1}^n  \ell_{ij}(u), 
\end{align}
or equivalently 
 \begin{align}
 \begin{split}
Y_{ij} (n)  &= \left( Y_{ij}(n-1) + \ell_{ij}(n) \right)^+, \quad n \in \bN,\\
Y_{ij}(0) &= 0,
\end{split}
\end{align}
we recover the Matrix CuSum in  \cite{oskiper2002online}.





\subsection{The Adaptive Matrix CuSum}
In this work, we  propose a novel  modification of the Matrix CuSum that  is obtained by resetting each $Y_{ij}$  whenever $Y_i$ becomes $0$.  Specifically, for each $i,j \in [K]$ with  $ j \neq i$,  we propose selecting $W_{ij}$ in \eqref{family4} as
\begin{align}\label{proposed_procedure}
Y'_{ij} (n) &\equiv \max_{R_i(n) \leq t \leq n} \sum_{u=t+1}^n  \ell_{ij}(u), 
\end{align} 
where $R_i(n)$ is, at time $n$, the most recent time that $Y_i$ was at $0$, i.e., 
\begin{align} \label{R}
R_i(n) &\equiv  \max \left\lbrace 0 \leq  t \leq  n   : Y_{i}(t) = 0 \right\rbrace.
\end{align}

This scheme  admits  the same recursive structure as the Matrix CuSum. Indeed,  for every $n \in \bN$ and $i, j \in [K]$ with $i \neq j$  we have 
\begin{align} \label{proposed_recursion}
Y^{\prime}_{ij}(n)  &=  \left(
             Y^{\prime}_{ij}(n-1)  + \ell_{ij}(n)  \right)^+  \cdot \mathds{1}  \left( \{ Y_i(n) > 0\} \right) , 
\end{align}
where $ Y^{\prime}_{ij}(0) = 0$. The main reason for proposing it is that it suppresses an  excessive use of pre-change data for  deciding the  post-change  distribution that, 
 as we will see later,  in  some cases characterizes  the Matrix CuSum and  leads to very large conditional probabilities of false identification. Indeed,  for each $ i \in [K]$, $ R_i(n)$ is an  estimate,  based on the data up to time $n \in \bN$, of the  time at which the density changes from $f$ to $g_i$, which goes back to \cite{hinkley1970} and   has been used in the literature of sequential change detection for different purposes (see, e.g.,  \cite{yang2017quickest}).  Therefore, discarding the data up to and including time $R_i(n)$ in the evaluation of $W_{ij}(n)$, for  each   $ j \neq i$ and $n \in \bN$,  allows the estimate  of the post-change density to be  based mostly on data from after the change, no matter when it occurs.    In  Section \ref{sec:analysis} we conduct a theoretical analysis that supports the above claims, but,  first,  in  Section \ref{sec:suppressing}, we  provide an intuitive explanation   of why the Matrix CuSum, as well as  the Generalized CuSum, can systematically fail to correctly isolate the post-change distribution in certain setups.\\

\noindent \textbf{\underline{Remark:}} An alternative sequential change diagnosis scheme  can be obtained by setting  $M_i(n) = R_i(n)$ in  \eqref{window_limited}.  
 However,  the calculation of $ W_i(n)$  in this case would require a number of operations of the order of  $n - R_i(n) $. Specifically, the computational cost and memory requirement of the procedure would be  dictated by a random variable with a long tail, which may not be desirable in practice.   On the other hand, applying  the window $ R_i(n)$ to the Matrix CuSum statistics leads to a  recursive structure \eqref{proposed_recursion}, which is another motivation for  our proposal of this scheme. \\

\section{The Effect of Pre-Change Data}\label{sec:suppressing}
In this section we explain intuitively how the Matrix CuSum and Generalized CuSum can systematically fail, in certain setups,  to correctly isolate the  post-change regime. To this end, we focus on the case that the true post-change density is $g_K$ and  there is another density, say $g_1$, such that $g_1$ and $g_K$ are closer together than $g_K$ and $f$, in the sense that
\begin{equation}
\sD (g_K \, || \, g_1) < \sD(g_K \, || \, f)
\end{equation}
or, equivalently,
\begin{align}
\Exp_K[\ell_{1}(n)] &>0 \quad  \forall \; n \in \bN.  \label{cond2}
\end{align}
This is the case, for example, when there is a change in the mean of a sequence of Gaussian random variables from $0$ to some positive number, i.e. 
\begin{equation} \label{Gaussian}
f = \cN(0,1) \quad \text{and} \quad  g_i = \cN(\theta_i, 1), \quad   i \in [K], \quad \quad 
\text{where} \quad 0 <\theta_1 < \ldots < \theta_K,
\end{equation}
 or in  the multichannel problem of Section \ref{subsec:multichannel} when more than one channel may be affected by the change. Consider, for simplicity, the case with $ d= 2$ channels, where $ K = 3$ and $ g_1, g_2, g_3$ are given by \eqref{multichannel_simultaneous}, i.e., $ g_i$ is the density when the change affects only channel $i$, where $ i \in \lbrace 1, 2 \rbrace$, and $ g_3 $ is the density when the change affects both channels.
Then, \eqref{cond2} holds, since 
\begin{align}
\begin{split}
\sD(g_3 \, || \, f ) &= \sD(q_2 \, || \, p_2) +  \sD(q_1 \, || \, p_1) \\
&> \sD(q_2 \, || \, p_2) = \sD(g_3 \, || \, g_1 ).
\end{split}
\end{align}

By the standard properties of the  CuSum statistics  and  condition \eqref{cond2}  it follows that  $Y_K$ will  tend to be close to 0 before the change and will increase after the change, and  that  this will also be the case for $Y_1$, although its growth after the change will be smaller than  that of  $Y_K$. As a result,  if also   $W_1$  is large, especially if it is  larger than  $W_K$, for a certain time period after the change,   it  becomes likely to  incorrectly identify  $g_1$ as the post-change density. This is the case, for example,  if, in addition to \eqref{cond2},  $g_1$ is the unique density that is the closest to $f$ before the change in the sense that 
\begin{equation}
 \sD(f \, || \, g_1) < \sD(f \, || \, g_j)   \text{ for all } j \in \{ 2, \ldots, K\},
\end{equation}
or equivalently 
\begin{align}
\Exp_\infty[\ell_{1j}(n)] &> 0 \text{ for all } j \in \{ 2, \ldots, K\}, \quad  n \in \bN.\label{cond1}
\end{align}
The latter condition is  satisfied for example in the  Gaussian mean shift problem \eqref{Gaussian}.  

To explain how condition \eqref{cond1} implies a large value of $W_1$, much larger than that of $W_K$,  for a certain period of time after the change,  for simplicity we focus on the simplest case  of  $ K = 2$ post-change alternatives and we  illustrate our explanation  
(Fig. \ref{fig:matrix_cusum_vs_proposed} \& \ref{fig:WLGC_paths})  in the context of  Gaussian mean shift problem. Thus,  in the following discussion, the  true post-change density is $g_2$.

\subsection{Matrix CuSum and the proposed  adaptive modification}  We start  with the Matrix CuSum, i.e., when   
\begin{equation}
W_1= W_{12}=Y_{12} \quad \text{and} \quad W_2= W_{21}=Y_{21}.
\end{equation}
Condition \eqref{cond1} then implies that,  before the change, $ W_2$ will tend to be close to 0 and  $W_{1}$ will    behave like a random walk with positive drift. As a result, the expected  value of $W_1$  at the time of the change  will be proportional to the change-point.  Therefore, even though $W_1$ will start decreasing after the change,  there will be  some  period  after the change that it will be large and, in fact, much larger than $ W_2 $, as the latter will start growing only after the change occurs. Most importantly, the length of this period  is increasing in the change-point. The longer it takes for the change to happen, the larger the  expected value  of $W_1$ at the time of the change, the longer it take for its values after the change to become small.  Thus, it is not possible to control this problematic behavior unless an upper bound  is imposed on the change-point.  
On the contrary,  using the  proposed adaptive window,   i.e., setting 
\begin{equation}
 W_1= W_{12}=Y'_{12} \quad \text{and} \quad W_2= W_{21}=Y'_{21},
\end{equation}  
effectively removes pre-change data from the decision-making process and prevents the erroneous growth of $W_1$ before the change, without affecting the correct behavior of $W_2$. 

\begin{figure}[tb]
 \begin{center}
\includegraphics[scale=1]{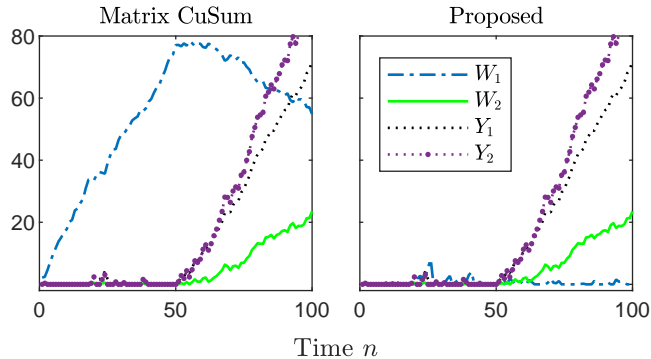}  
\caption{Sample paths of the Matrix CuSum and Adaptive Matrix CuSum statistics are plotted in the Gaussian mean shift problem where $ \nu = 50.$}
\label{fig:matrix_cusum_vs_proposed} 
\end{center}
\end{figure}

These points  are illustrated in  Fig. \ref{fig:matrix_cusum_vs_proposed},  where we see that 
$ Y'_{12}$, unlike $Y_{12}$, does not  grows before the change, whereas $Y'_{21} $ behaves very similarly to  $Y_{21}$.
\\

\subsection{ Generalized CuSum and window-limited modifications}  We continue with the Generalized CuSum,  i.e., when 
\begin{align} \label{gen_cum_stat}
\begin{split}
W_1(n) &= \max_{0 \leq k \leq n} \left\lbrace  \sum_{u = k+1}^n \ell_1(u) \wedge \sum_{u = k+1}^n \ell_{12}(u) \right\rbrace, \\
W_2(n) &= \max_{0 \leq k \leq n} \left\lbrace  \sum_{u = k+1}^n \ell_2(u) \wedge \sum_{u = k+1}^n \ell_{21}(u) \right\rbrace. 
\end{split}
\end{align}

\begin{figure}[tb]
\begin{center}
\includegraphics[scale=1]{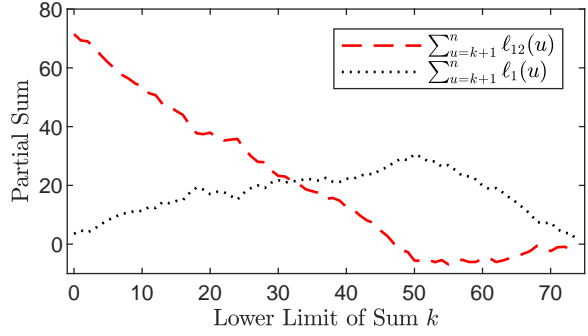}  
\caption{The partial sums $ \sum_{u = k+1}^n \ell_1 (u ), \sum_{u = k+1}^n \ell_{12}(u)$ are plotted as a function of the lower limit $k$ when $n = 75$ and $ \nu = 50$, in the Gaussian mean shift problem.} 
\label{fig:partial_sums} 
\end{center}
\end{figure}

Then,  standard properties of likelihood ratios and condition \eqref{cond1} imply that $W_2 $ does not grow before the change and starts increasing only after the change. On the other hand,  while $W_1$ does not grow before the change, as in the case of the Matrix CuSum, it does grow for a certain period of time \textit{after} the change occurs.  To see this, observe that for ${ k \geq \nu }$,  the sum $\sum_{u = k+1}^n \ell_{12}(u)$ is  generally small/negative, which is conducive to correct identification, while $ \sum_{u = k+1}^n \ell_1 (u )$ is generally positive. However, both of them can be large for $ k < \nu $, i.e., when pre-change data are used, as it is    illustrated in 
 Fig. \ref{fig:partial_sums}, where  these two sums are plotted  as  functions of $k$ when  $n=75$ and $\nu=50$. 
In this Figure we see that the first sum is maximized when $ k \approx \nu$, as expected, and the second is  negative when $ k \geq \nu $, as desired. However, both of them  are positive and quite large for $ k$ around  $35$, which  results in a  rather large value for  $W_1$ at time $n$.

\begin{figure}[tb]
 \begin{center}
\includegraphics[scale=1]{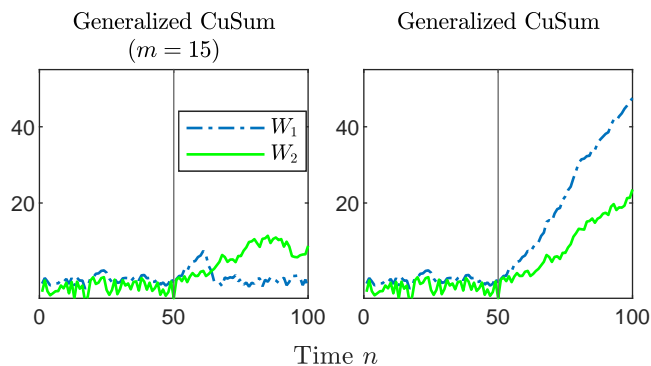}  
\caption{Sample paths of the Generalized CuSum statistics and the window-limited modification are plotted in the Gaussian mean shift problem where $ \nu = 50.$}
\label{fig:WLGC_paths} 
\end{center}
\end{figure}

To sum up, when applying the Generalized CuSum,   $W_1$  does not grow  before the change, but it does grow 
 for a certain time-period \textit{after the change}  during which it is  likely to be much larger  than $ W_2 $, as the latter grows only from post-change data. Most importantly,  the  length  of this problematic time-period increases with  $\nu$,  i.e., a later change-point results in  an  erroneous post-change growth of $W_1$  for a longer period of time, and therefore $W_1 $ becomes more likely to be larger than $W_2$. This point is illustrated in the right-hand side of Fig. \ref{fig:WLGC_paths}. 
 
This phenomenon is mitigated when using the window-limited Generalized CuSum, i.e., when restricting the maxima in \eqref{gen_cum_stat} to $n-m\leq k \leq n$,  where $m$ is some deterministic constant. Then, the   post-change growth of $W_1$ can last  for  at most $m$ observations,   after which  pre-change data  are removed from the calculation of the statistic. This is illustrated in the left-hand plot in Fig. \ref{fig:WLGC_paths} when $\nu=50$ and $m=15$.  On the other hand, as it shown in \cite{lai2000sequential},  $m$ should be selected large enough for the detection of the change to occur within the window, i.e., within $m$ observations after the change occurs. Therefore, the selection of $m$ is characterized by a trade-off, which is not clear how to resolve in a satisfactory way.   As we discussed in the final remark of the previous section, a solution to this problem would be to replace the fixed window $m$ by the adaptive window in \eqref{R}, i.e., to restrict the first maximum in \eqref{gen_cum_stat} to $R_1(n)\leq k \leq n$ and the second to $R_2(n)\leq k \leq n$. \\  


\noindent \textbf{\underline{Remark:}}  Conditions \eqref{cond2} and  \eqref{cond1} provide only a particular class of examples in which the Matrix CuSum and the Generalized CuSum can have difficulties in identifying the post-change. For example, condition   \eqref{cond1}  is not satisfied with strict inequality  in the case of the multichannel problem with  simultaneous faults. Nevertheless,   a similar reasoning  can be applied to argue that when  the change does not occur from the outset and it affects more than one channels, then these  two schemes are likely to make a wrong identification. In fact, this  is shown for the Matrix CuSum, and to a lesser extent the Generalized CuSum, in our simulation studies in Section \ref{sec:simulation}.

\section{Performance analysis}\label{sec:analysis}
In this section we establish some general results regarding the performance of sequential change diagnosis procedures of the form  \eqref{general_stopping_time}-\eqref{family}, where also each  statistic $W_i$ is of the form  \eqref{family4}. We focus, in particular,  on the proposed scheme, for which we establish (i)~an  upper bound on its Lorden delay,
(ii)~a novel non-asymptotic lower bound for its expected time to false alarm, (iii)~an upper bound for its  conditional probability of false isolation for a certain family of    change-points, and (iv)~an asymptotic optimality property as $\alpha$ goes to 0 sufficiently faster than $\beta$.

\subsection{Delay  analysis}
We start our delay analysis with 
 an upper bound for the expected delay of $\tau(b,h)$ when  $\nu=0$, i.e., when the change   has occurred once the  monitoring begins. For the following lemma, we introduce $I_i^*$ as the  minimum of the Kullback-Leibler information numbers in ${ \{I_{ij}:  j \in [K], j \neq i\} }$, i.e.,
 \begin{equation}
  I_i^* \equiv  \min_{j \in [K]: j \neq i}  I_{ij}.
 \end{equation}

\begin{lemma} \label{lem:new_lemma_for delay}
Fix $i \in [K]$ and suppose that,  for every $j \in [K]$ such that $j \neq i$, 
\begin{align} \label{summable2}
\sum_{n=1}^\infty \Pro_i \left( W_{ij}(n) \leq \rho n  \right) &< \infty \quad \text{for all} \; \rho < I_{ij}.
\end{align}
Then, for all $\delta > 0$, there is a constant $C_\delta>0$ such that
\begin{align} \label{delay under 0}
\Exp_i \left[ \tau_i(b,h) \right] \leq  \max \lbrace b / I_{i}, h / I^*_i \rbrace (1+\delta) +C_\delta
\end{align}
for all $b,h \geq 0.$
\end{lemma}

\begin{IEEEproof} Appendix C.\\
\end{IEEEproof}

Condition \eqref{summable2} is satisfied when $W_{ij}$ is either $Y_{ij}$ or $Y_{i}-Y_j$, which means that the upper bound \eqref{delay under 0} applies to both the Vector CuSum and the Matrix CuSum. Since, by definition, $Y'_{ij} \leq Y_{ij}$, it is not obvious that  \eqref{summable2} holds for  the proposed scheme, i.e., when $W_{ij}=Y'_{ij}$. In the following lemma we show that  this is indeed the case under mild moment conditions.

\begin{lemma}\label{lem:det_delay_bound} 
Fix  $i, j \in [K], j \neq i$. If  there is a   $ p > 1$  such that
\begin{equation}
\Exp_i[|\ell_{ij}(1)|^{p}] < \infty,
\end{equation}
 and an  $\epsilon>0$  so that 
 \begin{equation}
  \Exp_i \left[ \left(  \ell_i^-(1)  \right)^{2 + p/(p-1) + \epsilon} \right] < \infty, 
 \end{equation}
  
then \eqref{summable2} holds when $W_{ij}=Y'_{ij}$. 
\end{lemma}

\begin{IEEEproof}
  Appendix C.\\
\end{IEEEproof}

In the case  of the Matrix CuSum, i.e., when $W_{ij}=Y_{ij}$,  it is well known that the  Lorden delay of  $\tau(b,h)$ agrees with the expected delay of $\tau(b,h)$ when  $\nu=0$, i.e., for all $ i \in [K]$ and $b,h>0$,\begin{equation}
\cJ_i[\tau(b,h)]=\Exp_i[\tau(b,h)].
\end{equation}
In the following lemma we show that  this is also the case for the proposed scheme, i.e., when $W_{ij}=Y'_{ij}$, despite the presence of the adaptive window. We emphasize  that this is \textit{not} the case for the  Vector CuSum, whose   worst-case delay analysis  was   conducted  in \cite{nikiforov2000simple,nikiforov2003lower} using  Pollak's criterion,  not Lorden's.

\begin{lemma}\label{lem:worst_case}
If $W_{ij}=Y'_{ij}$ for every $i,j \in [K], j \neq i$, then 
\begin{align} \label{worst case at 0}
\cJ_i[\tau(b,h)] = \Exp_i[\tau(b,h)] \quad \text{ for all }   i \in  [K], \; b, h\geq 0.
\end{align}
\end{lemma}

\begin{IEEEproof}
Fix  $ i \in  [K]$ and  $b, h\geq 0$.  By the definition of $\tau(b,h)$ and  recursions \eqref{cusum_recursion} and \eqref{proposed_recursion} it follows  that, for any $\nu \in \bN_0$, the conditional expectation 
\begin{equation}\label{eqn:conditional_expectation}
\Exp_{\nu, i} \left[\tau(b,h) - \nu \; | \; \mathcal{F}_{\nu }, \tau(b,h) > \nu  \right]
\end{equation}
depends on $\cF_\nu$ only through
\begin{equation}
  \cY(\nu) \equiv \{ Y_k(\nu), Y_{kj}'(\nu), \; k, j \in [K], j \neq k \}.
\end{equation}
  Thus, for every $\nu \in  \bN_0$, 
\begin{align}
\begin{split}
& \Exp_{\nu, i} \left[ \tau(b,h) - \nu \; | \; \cF_{\nu}, \tau(b,h)  > \nu \right] \\
&= \Exp_{\nu, i} \left[\tau(b,h)  - \nu  \;|\; \cY(\nu), \tau(b,h)  > \nu \right].
\end{split}
\end{align} 
From recursion \eqref{proposed_recursion}   it also follows  that the latter conditional expectation is decreasing, thus, it is maximized when all components of  $\cY(\nu)$ are as small as possible, i.e., $0$. When this happens,  this conditional expectation becomes equal to $\Exp_i[\tau(b,h)]$, since   $ X_{\nu + 1}, X_{\nu + 2}, \ldots$ are independent of $\cF_\nu$ and have the same distribution under~$ \Pro_{\nu,i}$ as $X_1, X_2, \ldots$ under~$\Pro_i$. \\


\end{IEEEproof}

Combining the above results, we establish an upper bound for the Lorden delay of the proposed scheme. 

\begin{theorem}\label{thm:delay_bound}
If $W_{ij}=Y'_{ij}$ for every $i,j \in [K], j \neq i$ and the conditions of Lemma \ref{lem:det_delay_bound} hold, then  for every $\delta>0$ there is a  constant $C_\delta>0$ so that, for every  $b>0$ and $h>0$, 
\begin{align} \label{Lorden delay}
\cJ_i \left[ \tau(b,h) \right] \leq  \max \lbrace b / I_{i}, h / I^*_i \rbrace (1+\delta) +C_\delta.
\end{align}
\end{theorem}

\begin{IEEEproof}
The proof follows directly by combining the lemmas of this subsection.\\
\end{IEEEproof}

\subsection{False alarm control}

 Independently of the  choice of the statistics $W_{i}$, $i \in [K]$, by the definition of $\tau(b,h)$ it is clear that, for any  $b, h\geq 0$, it holds that $\tau(b,h) \geq \sigma (b)$ and, in view of  \eqref{CUSUM_ARL}, 
\begin{align*}
\Exp_\infty \left[ \tau (b,h)\right] &\geq   e^{b}/K.
\end{align*} 
 Thus, it is always possible to guarantee that the false alarm rate of $\tau (b,h)$ does not exceed an arbitrary $\alpha \in (0,1)$ by simply selecting  $b$  to be larger than or equal to $b_\alpha$, which was defined in \eqref{eqn:b_alpha}.  Indeed,  this is the approach typically taken in the literature, at least for the purpose of  analysis.   However, as we show in the next lemma,  if the $W_{ij}$s are stochastically small under $\Pro_\infty$, then we also have  a lower bound for the expected time to false alarm that is  independent  of  $b$ and exponential in $h$. 

\begin{lemma}\label{lem:false_alarm_bound}
If, for every  $ i, j  \in [K]$ such that $j \neq i$,    there are  $Q_{ij}, q_{ij} > 0$  so that,   for all $ n \in \bN$, 
\begin{equation} \label{tail}
\Pro_{\infty} \left(  W_{ij}(n) \geq x \right)  \leq Q_{ij}\,  e^{-q_{ij} x}, \quad \forall  x \geq 0,
 \end{equation}
 then, for any $b>0$ and $h>0$, 
\begin{align}
\Exp_\infty \left[ \tau(b,h) \right]  &\geq \frac{1}{2}  \left( \sum_{i=1}^K \min_{j \neq i}  \left\{ Q_{ij} e^{-q_{ij} h} \right\} \right)^{-1}.
\end{align} 

\end{lemma}

\begin{IEEEproof}
Appendix D.\\

\end{IEEEproof}

Condition \eqref{tail}  is  satisfied for all $n \in \bN$ in the case of the Vector CuSum, i.e., when $W_{ij}=Y_i-Y_j$. Indeed, for  any such $i,j$,   $n \in \bN$ and $x \geq 0$,  
\begin{equation}
\Pro_{\infty}(Y_i(n) - Y_j(n) \geq x) \leq \Pro_{\infty}(Y_i(n)  \geq x) \leq e^{-x},
\end{equation} 
where the second inequality is a well-known property of the CuSum statistic (restated in  Appendix A).   On the other hand,  condition \eqref{tail} is not, in general, satisfied in the case of the Matrix CuSum, i.e., when    $W_{ij}=Y_{ij}$.  Indeed, when $\Exp_\infty[\ell_{ij}(n)] >0$,  $Y_{ij}$  essentially behaves as a random walk with positive drift before the change.  In the next theorem we show that  condition \eqref{tail} is always satisfied for the proposed scheme as long as condition \eqref{assum:A1_prime} holds.

\begin{theorem}\label{thm:uniform_exponential_bound}
 If condition \eqref{assum:A1_prime} holds, then  for every $i, j \in [K]$ such that $j \neq i$ there are $Q_{ij}, q_{ij} >0$, such that  \eqref{tail} holds for all $ n \in \bN_0$  when $W_{ij}=Y'_{ij}$.

\end{theorem}

\begin{IEEEproof}
Appendix E.\\
\end{IEEEproof}

In summary, the results of this section reveal  that   both the proposed diagnosis scheme and  the Vector CuSum can control the false alarm rate below any user-specified value $\alpha$ with a suitable selection of  threshold $h$ alone, independently of the value of $b$.

\subsection{False isolation control}
To establish our main result regarding false isolation control, we need  $W_{ij}$ to be   pathwise smaller than or equal to $Y_{ij}$, i.e., 
\begin{equation} \label{pathwise bound}
W_{ij}(n) \leq Y_{ij}(n) \quad \text{for all} \quad 
n \in \bN
\end{equation} 
for all $i, j \in [K], j \neq i$. This holds, trivially, for the Matrix CuSum and the proposed scheme,  whereas it also holds for the Vector CuSum.
 In fact, in our analysis, we use  the same technique as in the  proof for the false isolation control of the Vector CuSum in  \cite[Theorem 2]{nikiforov2000simple}, but, as we will see, our conclusions are somewhat different. 
 
To state our results in the appropriate generality, we need to introduce, 
for each $i, j \in [K], j \neq i$, the fictitious statistic that  utilizes the CuSum statistic $Y_{ij}$, instead of $W_{ij}$, after $\nu$, i.e., 
\begin{align} \label{W_nu}
  W_{ij}^{\nu}(n) \equiv 
        \begin{cases}
            W_{ij}(n), & \quad   n \leq \nu, \\
           \left(W^{\nu}_{ij}(n - 1) + \ell_{ij}(n)  \right)^+, & \quad  n > \nu,
        \end{cases}
 \end{align}
and the stopping time
\begin{align}\label{tau_nu}
 \tau^\nu_{ij}(h)  &\equiv \inf \lbrace n > \nu : W^\nu_{ij}(n) \geq h \rbrace.
 \end{align}
Specifically, the statistic $ W_{ij}^{\nu}(n)$ is equal to $ W_{ij}(n)$ before the change occurs, and behaves like a CuSum statistic with initialization $ W_{ij}(\nu)$ after the change-point, $\nu$.

\begin{theorem}
\label{th:family_misspecification_bound_theorem}
Fix  $b >0 $, $h \geq 1$,  $\nu \in \bN_0$,
$i, j \in [K]$ with $i \neq j$. Let $ C_{ij} > 0$, $ c_{ij} \in (0,1)$. If  $W_{ij}$ satisfies \eqref{pathwise bound} and 
\begin{equation} \label{condition on the change-point} 
\Pro_{\infty}(W_{ij}(\nu) \geq x \, | \, \tau(b,h) > \nu) \leq C_{ij}e^{-c_{ij}x}  \qquad  \text{for all } \;  x \in (0, h],
\end{equation}
and 
\begin{align} 
  \Pro_{\nu, j} \left( \tau_{ij}^{\nu}(h) > \tau_j(b,h) \, |\,  \tau(b,h) > \nu \right) > 0, \label{condition technical} 
  \end{align}
then there exists a function   
$$ \phi_{ij}: (0, \infty) \to (0, \infty)$$
   such that  $\phi_{ij}(h) \to 0$ as $h \to \infty$,  and 
\begin{align} \label{show}
\begin{split}
 \Pro_{\nu, j} (\widehat{\tau}(b,h) & =i  \,| \, \tau(b,h) > \nu )   \\
 &  \leq   e^{-h}  \, \Exp_{\nu, j}[\tau_j(b,h) - \nu \,|\, \tau(b,h) > \nu]  + \frac{C_{ij}}{1 - c_{ij}}e^{-c_{ij}h} \left( 1 + \phi_{ij}(h) \right).
\end{split}
 \end{align}

   \end{theorem}

\begin{IEEEproof}
Appendix F. \\
\end{IEEEproof} 

\noindent \underline{\textbf{Remark}}:  If condition \eqref{condition on the change-point}  holds with $c_{ij} \geq 1 $, then  a tighter upper bound can be derived with minor adjustments to the proof.\\

To explore the implications of  Theorem   \ref{th:family_misspecification_bound_theorem} we need to introduce some additional notation. Thus,  for any family of real numbers
\begin{equation}
 \bs{C} \equiv \left\lbrace c_{ij} \in (0,1), C_{ij}>0: \, i, j \in [K], i \neq j \right\rbrace,
\end{equation}
we introduce  the set of numbers in $\bN_0$ for which \eqref{condition on the change-point} holds for all $ b > 0, h \geq 1$:
\begin{align}
N_{\bs{C}} \equiv \big\lbrace n \in \bN_0 : \Pro_{\infty}(W_{ij}(n) \geq x \, | \, \tau(b,h) > n) \leq C_{ij}e^{-c_{ij}x} \quad  \forall \,  x \in (0, h], \, b > 0, \, h \geq 1 \big\rbrace .
\end{align}

\begin{lemma}
If  $W_{ij}(0)= 0$  for every $i \neq j$,  then $0\in N_{\bs{C}}$ for any $\bs{C}$.  If $m \in \bN$ and condition \eqref{tail} holds for every $n \in [m]$, then  there is a $ \bs{C}$ so that  $[m] \subseteq N_{\bs{C}}$.
\end{lemma}

\begin{IEEEproof}
The first claim is obvious. For the second, observe that, for every $n \in \bN$, 
\begin{equation}
\Pro_{\infty}(W_{ij}(n) \geq x \, | \, \tau(b,h) > n) \leq \frac{\Pro_{\infty}(W_{ij}(n) \geq x )}{ \Pro_{\infty}( \tau(b,h) > n)}.
\end{equation}
By standard properties of the CuSum  statistics,
\begin{equation}
\Pro_{\infty}( \tau(b,h) > n) \geq  \Pro_{\infty}( \sigma(b) > n) > 0.
\end{equation}
Therefore,  if condition \eqref{tail} holds for every $n \in [m]$, then there is a $ \bs{C}$, which depends on $m$, such that $[m] \subseteq N_{\bs{C}}$. \\
\end{IEEEproof}

\noindent \underline{\textbf{Remark}}: As we showed  earlier,
 \eqref{tail} holds \textit{uniformly in $n$} for  the proposed procedure and   the Vector CuSum. In view of this, we conjecture that, for these two schemes,   there is a $ \bs{C} $ such that $ N_{\bs{C}} = \bN_0$.   However, this cannot be shown with the previous argument, since 
 \begin{equation}
  \inf_{n \in \bN_0} \Pro_{\infty}( \sigma(b) > n) =
 0. 
 \end{equation}

We next specialize Theorem   \ref{th:family_misspecification_bound_theorem} to the proposed scheme.

\begin{corollary}\label{coro: h_choice}
Let  $W_{ij}=Y'_{ij}$ for all $i, j \in [K]$ such that $ j \neq i$. Suppose that condition \eqref{assum:A1_prime} holds, and the  conditions of Lemma \ref{lem:det_delay_bound} hold  for all $i, j \in [K]$ such that $ j \neq i$. Then, for any  $\alpha, \beta \in (0,1)$ and any $ \bs{C}$, we can select $h$ large enough so that 
\begin{equation}
(\tau (b_\alpha, h),  \widehat{\tau}(b_{\alpha}, h)) \in \cC(\alpha, \beta, N_{\bs{C}}).
\end{equation}
 \end{corollary}

\begin{IEEEproof}
Appendix F. \\
\end{IEEEproof}






\noindent \underline{\textbf{Remark}}: A similar corollary of Theorem \ref{th:family_misspecification_bound_theorem} can be stated for  the Vector CuSum, i.e.,  when  ${W_{ij}=Y_i-Y_j}$ for every $i \neq j$. Indeed, the conditional expected delay in the upper bound of \eqref{show} can be upper bounded using \cite[Theorem 1]{nikiforov2000simple}. However, it is not clear whether condition \eqref{condition technical} is satisfied for all $ \nu \in \bN_0$, at least for $h$ large enough, in the case of the Vector CuSum. Even if it is, the corresponding corollary would not guarantee false isolation control uniformly in $ \nu \in \bN_0 $, as it is  implied by \cite[Theorem 2]{nikiforov2000simple}.  The reason is that in the proof presented in \cite{nikiforov2000simple} the conditioning in \eqref{condition on the change-point} is ignored, but it is not clear to us whether this is indeed possible.

\subsection{An asymptotic optimality property}

We close this section by showing that  the proposed procedure  achieves the infimum in \eqref{infimum} to a first-order asymptotic approximation as $\alpha$ goes to  $ 0 $ sufficiently faster than $\beta$, for a given $ \bs{C}$.

\begin{theorem}\label{th:optimality}
Let $W_{ij}=Y'_{ij}$ for all $i, j \in [K]$ such that $ j \neq i$ and suppose that \eqref{assum:A1_prime}, the conditions of Lemma \ref{lem:det_delay_bound} hold.
Then, for any $ \bs{C}$, there is family of thresholds $(h_\beta)$ such that, for all $ i \in [K]$, 
\begin{align}
& \sup_{\nu \in N_{\bs{C}}} \Pro_{\nu, i}(\widehat{\tau}(b_{\alpha}, h_\beta) \neq  i| \tau(b_{\alpha}, h_\beta) > \nu)  \lesssim \beta, \quad \text{and}\\ 
& \cJ_i[\tau(b_{\alpha}, h_{ \beta})]  \sim \inf_{(T, D) \in \cC(\alpha, \beta, N_{\bs{C}})} \cJ_i[T], 
\end{align}
 as $ \alpha, \beta \to 0 $ such that   \begin{equation}
  \log | \log \beta |  \ll \log | \log \alpha |  \ll| \log \beta |.
\end{equation} 

\end{theorem}

\begin{IEEEproof}
Appendix F. \\
\end{IEEEproof}

\section{Design } \label{sec:design}
The  formulation of Section \ref{sec:state} suggests
that  thresholds $b$ and $h$ of a sequential change diagnosis procedure of the form  \eqref{general_stopping_time}-\eqref{family}  should be selected so that its false alarm rate does not exceed  $\alpha$,  i.e., 
\begin{equation}\label{eqn:false_alarm_rate}
\Exp_\infty[\tau(b,h)] \geq 1/\alpha,
\end{equation}
where  $\alpha$ is a  user-specified number in $(0,1)$, and its worst-case conditional probability of false isolation   does not exceed $\beta$, i.e., 
\begin{equation} \label{optimizers}
 \max_{j \in [K]} \sup_{\nu \in \bN_0 }\Pro_{\nu, j}(\widehat{\tau}(b,h) \neq j | \tau(b,h) > \nu) \leq \beta,
\end{equation}
where   $\beta$ is a user-specified number in $(0,1)$.  As we have argued earlier, this  quantity cannot, in general, become arbitrarily small  for each scheme we consider in this work, even for very large  values of  $b$ and $h$. Moreover, even for the schemes for which it  can be controlled, the change-point for which the  supremum  is achieved is unknown. Therefore, in order to satisfy the second constraint, if it is even possible to do so, the  conditional probability
$$\Pro_{\nu, j}(\widehat{\tau}(b,h) \neq j | \tau(b,h) > \nu)
$$
needs to be estimated, for  various values of  $(b,h)$, for a sufficiently large number of  $\nu$s.

Motivated by these  considerations,  in this section we propose a  novel method for selecting the thresholds $b$ and $h$, which is more  computationally efficient than the one described above, as it does not involve computation of  conditional probabilities of false isolation. For this, we propose  selecting $b$ and $h$ to satisfy, in addition to the false alarm constraint \eqref{eqn:false_alarm_rate}, 
\begin{equation}\label{eqn:delay_constraint}
\cJ_i[\tau(b,h)] \leq r  \, \max_{j \in [K]} \mathcal{L}_j(\alpha), \quad \text{ for all } i \in [K],
\end{equation}
for some  user-specified  constant $r>1$,  where 
$ \cL_j(\alpha)$ is  the optimal worst-case Lorden delay in $ \mathcal{C}(\alpha) $ when $g = g_j$,  i.e., 
\begin{align} \label{optimal_cusum_performance}
 \cL_j(\alpha) \equiv \inf_{(T,D)  \in \mathcal{C}(\alpha)} \mathcal{J}_j[T] =\Exp_{j}[\sigma_j(b_j(\alpha))], 
 \end{align}
$b_j(\alpha)$ being  the threshold $b$  for which  
$\Exp_{\infty}[\sigma_j(b)]=1/\alpha$ (\cite{moustakides1986optimal}). In other words, we propose selecting  $b$ and $h$ so that the maximum Lorden delay of $ \tau(b,h)$ can increase at most by a factor $r$  relative to that of the optimal procedure that knows a priori the true post-change distribution.

More formally, given   user-specified values of   $ \alpha \in (0,1)$ and  $r>1$,  the set  of thresholds $(b,h)$ for which the average time to false alarm and detection delay conditions hold  can be expressed as follows:
\begin{equation}
\mathcal{S}(\alpha, r)  \equiv \bigcap_{i \in [K]} \mathcal{D}_i(\alpha, r )\cap \mathcal{A}(\alpha),
\end{equation}
where $\mathcal{A}(\alpha)$ is the set of pairs for which the false alarm constraint is satisfied, i.e., 
\begin{equation}
\mathcal{A}(\alpha) \equiv \left\lbrace (b,h) : \Exp_{\infty}[\tau(b,h)]\geq  1 / \alpha \right\rbrace,
\end{equation}
and, for $i \in [K]$,  
$\mathcal{D}_i(\alpha, r)$  is the set of pairs satisfying \eqref{eqn:delay_constraint}, i.e., 
\begin{align}
 \mathcal{D}_i(\alpha, r) &\equiv \left\lbrace (b,h) : \cJ_i[\tau(b,h)] \leq r \max_{j \in [K]} \cL_j(\alpha)\right\rbrace.
\end{align}
Once  $\mathcal{S}(\alpha, r)$ is computed, a rather natural and intuitive selection for $b$ and $h$ in $\cS(\alpha, r)$  is   to \textit{select $h$ as large as possible within this region, and to subsequently select $b$ as large as possible given this choice of $h$.} While this selection is not guaranteed to  optimize \eqref{optimizers}, our simulations studies in the next section suggest that it can, at least, provide a  good approximation to the true optimizer, whose computation is not feasible.





\subsection{Computational considerations}
The implementation of the proposed design requires, first of all,   the computation, for every $i \in [K]$,  of  the optimal Lorden delay when the post-change is $g_i$, namely $\mathcal{L}_i(\alpha)$, defined in   \eqref{optimal_cusum_performance}. 
This task can be performed easily using  Monte Carlo simulation, or using existing approximations, such as Siegmund's corrected Brownian approximations (\cite{siegmund1979corrected}). 

Second, it  requires the computation of the  expected time to false alarm and the Lorden delay of the scheme of interest, $\Exp_\infty[\tau(b,h)]$ and $\cJ_{i}[\tau(b,h)]$, for every ${i \in [K]}$ and various values of $b, h$.  This task  is particularly simple for procedures that satisfy 
 \eqref{worst case at 0}, i.e., for which 
 \begin{equation}
 \cJ_{i}[\tau(b,h)]= \Exp_i[\tau(b,h)], \quad i \in [K],
 \end{equation}
such as the min-CuSum, the Matrix CuSum and the proposed scheme. Indeed, for these schemes, one simply needs  to simulate paths of the sequence $(X_n)$ under $\Pro_\infty$ and $\Pro_i$, $i \in [K]$,  which allow 
the estimation of $\Exp_\infty[\tau(b,h)]$ and $\Exp_i[\tau(b,h)]$, for ${i \in [K]},$  
simultaneously for a grid of values of $ (b,h)$. 

On the other hand, for procedures that do  not satisfy \eqref{worst case at 0}, such as  the Generalized CuSum, its window-limited version, and the Vector CuSum, the Lorden delay can  be impossible to estimate, and only a biased, low estimate of  the Pollak delay may be obtainable in practice. Thus, it is not in general possible to have a fair comparison  between a procedure that satisfies \eqref{worst case at 0} and a procedure that does not.   

In view of this, in  our simulation studies in the next section  we evaluate, for simplicity,  the expected delay  for these procedures only when $\nu=0$, keeping in mind that  this  can be much smaller than the actual Lorden delay for procedures that do not satisfy \eqref{worst case at 0}.


\section{Simulation studies}\label{sec:simulation}
In this section we compare the min-CuSum, the Matrix CuSum, the Vector CuSum, the window-limited Generalized CuSum, and the proposed procedure in the  multichannel setup of Subsection \ref{subsec:multichannel}. We focus on the case that there are two channels  $(d=2)$, with Gaussian pre-change and post-change distributions, specifically 
\begin{equation}\label{eqn:p_and_q}
 p_i \equiv \cN(0,1), \quad q_i \equiv \cN(1,1),
 \end{equation}
 and we consider two setups for this problem. In the first one the change can happen in only one of the two channels, i.e.,  \eqref{multichannel_single} holds with ${K = d = 2}$. In the second, the change can  also happen in both channels simultaneously, in which case $K=3$ and  \eqref{multichannel_simultaneous} holds. 

\subsection{Design} 
  In both setups, we design each scheme according  the method  described in the previous section. Specifically, we fix $ \alpha = 1 \%$  and  consider several values of $r$.  We consider $b$ in increments of $0.01$ starting from $0$ and $h$ in increments of $ 0.05$ starting from $ 0.05$. 
  
To implement each of the above  schemes, we first  estimate,  for each $i \in [K]$,  the optimal Lorden  delay in   \eqref{optimal_cusum_performance}, with  $b_i(\alpha)$ replaced by the smallest  $b$ in our grid  for which  the estimate of 
 $\Exp_\infty[\sigma_i (b)]$ exceeds $1/\alpha$. The estimated thresholds and detection delays are provided in Table \ref{table:traditional_multichannel_table} and Table \ref{table:multichannel_table}.

 Subsequently,  for each diagnosis procedure, we estimate  the regions $ \mathcal{D}_1(\alpha, r), \mathcal{D}_3(\alpha, r), $ and $\mathcal{A}(\alpha)$, thereby obtaining  $ \cS(\alpha, r)$ for various values of $r$.  
This process is illustrated with figures for  the proposed procedure and for the Matrix CuSum when $r=1.3$ and when $r=2$   in Appendix G.
 
For the Monte Carlo estimation of each of these quantities, i.e., $\Exp_{\infty}[\sigma_i (b)]$, $\Exp_{\infty}[\tau(b,h)]$, and $\Exp_i[\tau(b,h)], i \in [K],$  we simulate  $ {L=5 \times 10^4} $ paths of relevant statistics under $ \Pro_1 $ and $ \Pro_3$ and $ 0.1\cdot L$ paths under $ \Pro_{\infty}$. We do not need to simulate any paths under $ \Pro_2 $ due to its symmetry with $ \Pro_1$. The standard error of the estimate of $\Exp_\infty[\sigma_i (b)]$ is less than $ 1.4 \%$ in each case. 
  
\begin{table}[b]
\caption{Multichannel problem with single faults - optimal worst-case detection delay:\\ Estimate (Standard Error)}
\label{table:traditional_multichannel_table}
\begin{center}
\begin{tabular}{ cccc } 
  & $b_i(\alpha)$ &$ \cJ_{i}[\sigma_i(b_i(\alpha)]$ \\ 
  \hline
 $i =1$ & {2.85} &  6.0797 (0.0166)\\ 
 \hline
\end{tabular}
\end{center}
\end{table}

\begin{table}[bt]
\caption{Multichannel problem with simultaneous faults - optimal worst-case detection delay: Estimate (Standard Error)}
\label{table:multichannel_table}
\begin{center}
\begin{tabular}{ cccc } 
  & $b_i(\alpha)$ &$ \cJ_{i}[\sigma_i(b_i(\alpha))]$ \\
  \hline
 $i =1, 2$ & {2.85} & 6.0965 (0.0165)   \\ 
 $i = 3$ & 3.04 & 3.7450 (0.0097)   \\
 \hline
\end{tabular}
\end{center}
\end{table}

\subsection{Comparisons}
Since  it is not feasible to compute the worst-case conditional probability of misidentification in \eqref{optimizers}, we compare the various schemes based on their worst-case conditional probabilities of misidentification when  $ {\nu  \in \lbrace 0, 10, \ldots, 50 \rbrace}$, i.e., 
\begin{equation}\label{eqn:simulated_probabilities}
  \max_{j \in [K]} \max_{\nu \in \lbrace 0, 10, \ldots, 50 \rbrace}\Pro_{\nu, j}(\widehat{\tau}(b,h) \neq j | \tau(b,h) > \nu) 
\end{equation} 

For each scheme, we compute this quantity in two ways as far as it concerns the selection of $b$ and $h$: 
\begin{itemize}
\item  using the proposed method in the previous section, i.e., the largest value of $h$ in  $\mathcal{S}(\alpha, r)$,  and subsequently the largest value of  $b$ given this choice of $h$, 
\item using the values of $b$ and $h$ that optimize \eqref{eqn:simulated_probabilities} in $\cS(\alpha,r)$. 
\end{itemize}

The first method requires no additional simulation, after computing $ \cS(\alpha, r)$, whereas for the second task  we simulate $L$ paths of relevant statistics for the Monte Carlo estimation of each of the probabilities in \eqref{eqn:simulated_probabilities}.

The two different methods for selecting $b$ and $h$  provide indistinguishable results in all cases apart from that of the Matrix CuSum in the setup of a single fault,  where  \eqref{eqn:simulated_probabilities} is optimized by 
a smaller $h$ and a larger $b$ compared to the ones obtained by the method of the previous section.  Although the difference is not dramatic, we present the results with the second choice of thresholds,  making sure that all comparisons are fair for the  Matrix CuSum.

\begin{figure}%
    \centering
    \subfloat[\centering Without simultaneous faults]{{\includegraphics[width=6.8cm]{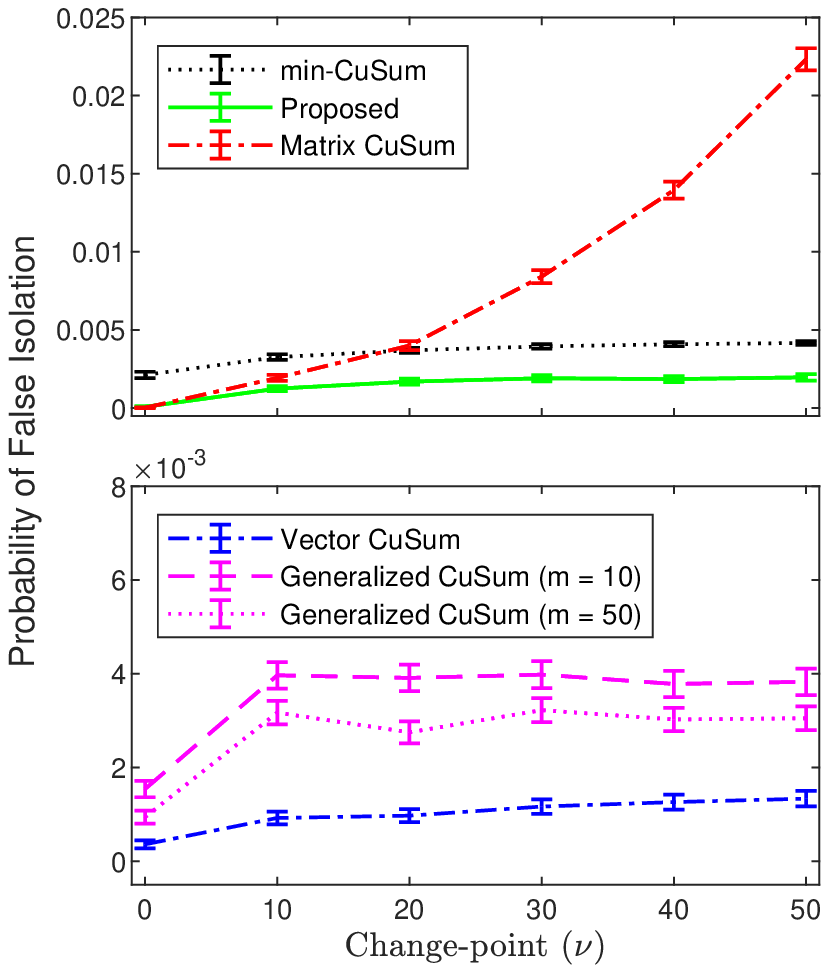} }}%
    \qquad
    \subfloat[\centering With simultaneous faults]{{\includegraphics[width=6.8cm]{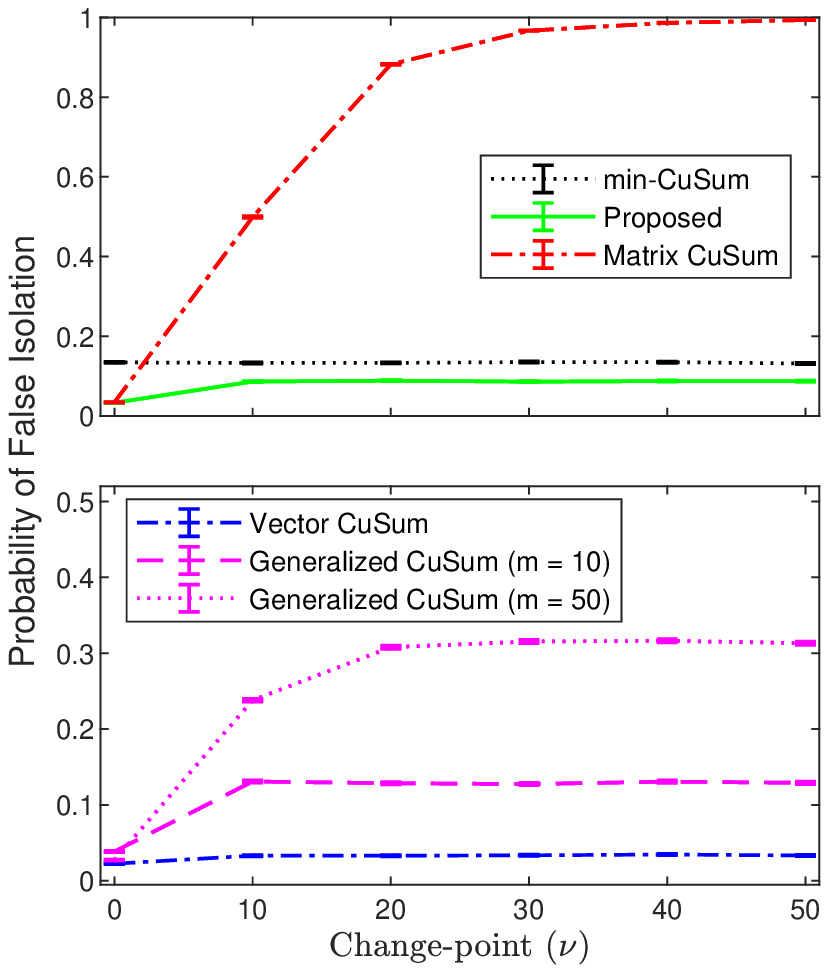} }}%
    \caption{Worst-case probabilities of false isolation (with respect to the type of change) plotted against the change-point $\nu$ with $ \alpha = 1 \% $ and $ r= 2 $ in the multichannel problem without and with simultaneous faults, respectively. Standard errors of probability estimates provided as error bars.}%
    \label{fig:stability}%
\end{figure}

In   Fig. \ref{fig:stability}  we plot, for each scheme,  the worst-case probability of false isolation 
$$  \max_{j \in [K]} \Pro_{\nu, j}(\widehat{\tau}(b,h) \neq j | \tau(b,h) > \nu) $$
against ${ \nu \in \lbrace 0, 10,\ldots, 50 \rbrace}$, when $r=2$,  for the single fault  and  the simultaneous fault setup, respectively. These graphs suggest that 
\eqref{eqn:simulated_probabilities} is a   reasonable proxy for  \eqref{optimizers}  for all schemes, apart from the Matrix CuSum. This means that using  \eqref{eqn:simulated_probabilities} instead of \eqref{optimizers}  in our comparisons is favorable to the Matrix CuSum.

\subsection{Results}
\begin{figure}%
    \centering
    \subfloat[\centering Without simultaneous faults]{{\includegraphics[width=6.8cm]{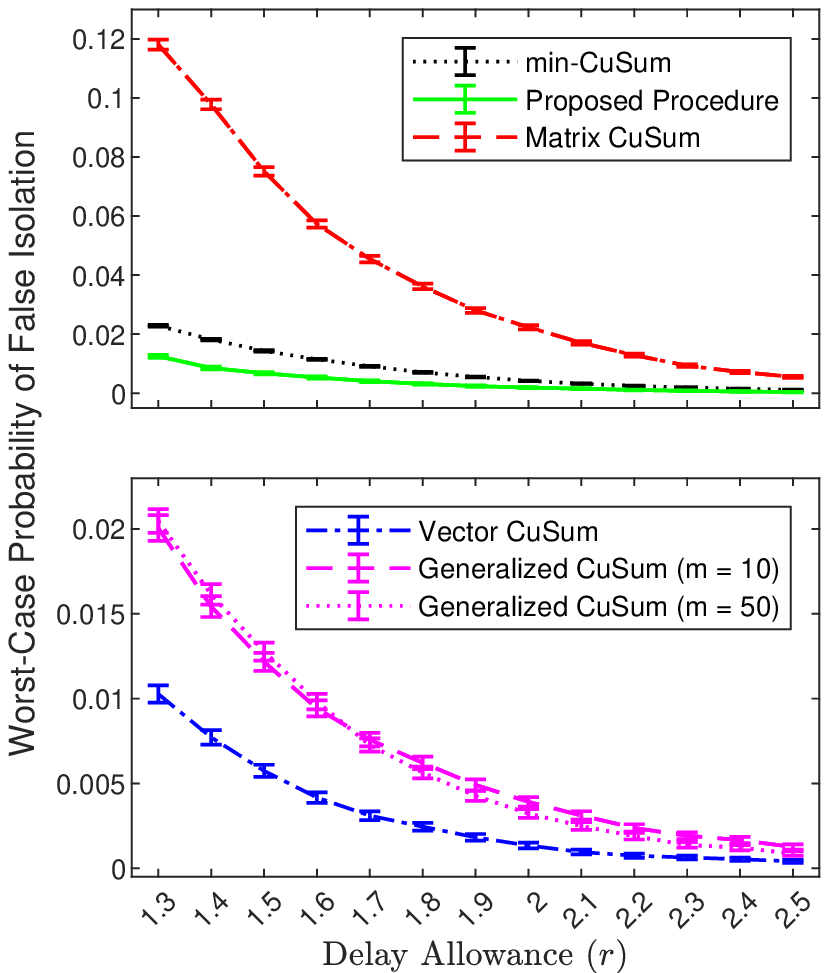} }}%
    \qquad
    \subfloat[\centering With simultaneous faults]{{\includegraphics[width=6.8cm]{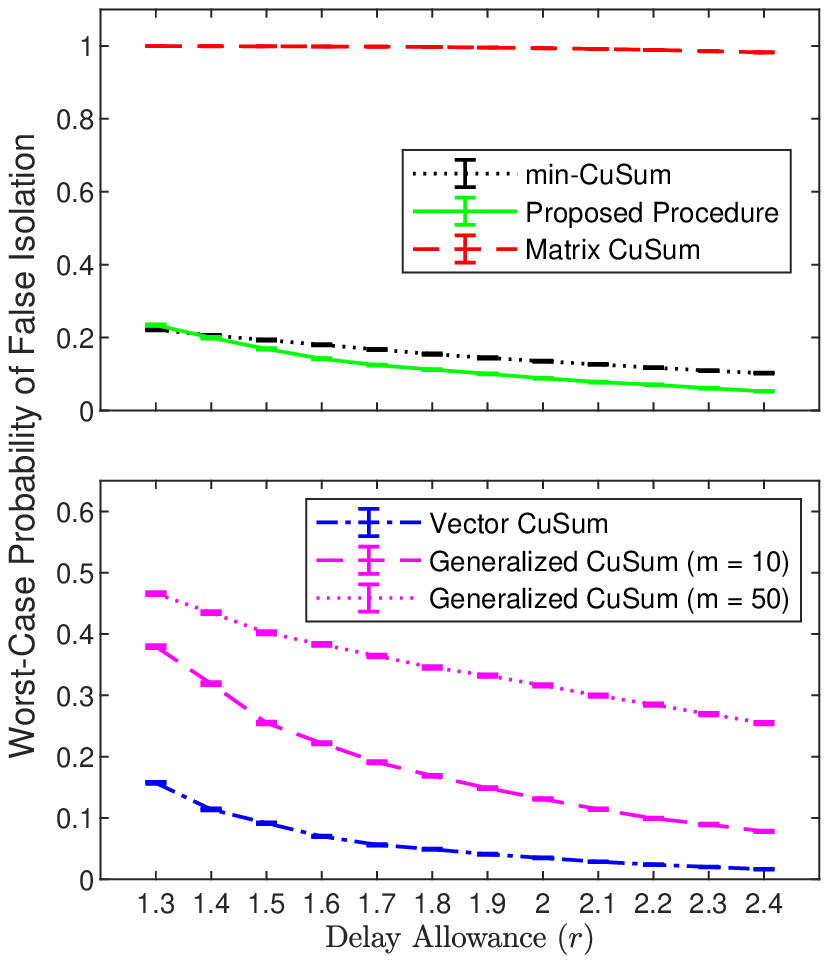} }}%
    \caption{Worst-case probabilities with respect to the change-point and post-change density as a function of the detection delay allowance factor $r$ in the multichannel problem without and with simultaneous faults, respectively. Standard errors of probability estimates provided as error bars.}%
    \label{fig:comparison}%
\end{figure}

In  Fig. \ref{fig:comparison}  we plot, for each scheme,   the  worst-case conditional probabilities of misidentification in  \eqref{eqn:simulated_probabilities}  against $r$ for the single fault and the simultaneous faults setup, respectively.  As explained  in the previous section, we compare separately the schemes that satisfy \eqref{worst case at 0}, i.e.,  the  proposed procedure, Matrix CuSum, and min-CuSum, and   those that do not, i.e., the Vector CuSum  and  the Generalized CuSum, since the constraints have different interpretations between the two groups.

From these graphs we observe that, in  both setups and for all values of $r$,    the proposed scheme performs as well or  better than the min-CuSum, and the latter performs better than the Matrix CuSum. However, in  the single fault setup  the differences between these  schemes vanish  as $r$ increases. On the other hand,  in the simultaneous fault setup the worst-case probability of the Matrix CuSum is very close to 1 even for large values of $r$, whereas for the two other schemes it is not much larger than $0.2$ even for  small values of $r$.

Moreover, we observe that, in  both setups and for all values of $r$,   the Vector CuSum performs  better than  the window-limited Generalized CuSum, independently of the choice of window size, $m$, in the latter.  However, in  the single fault setup  the choice of $m$ does not make any practical difference, whereas  the difference between the two  schemes vanishes  as $r$ increases. On the other hand,  in the simultaneous fault setup a larger $m$ leads to substantially worse misidentification probability for the Generalized CuSum and, even with a smaller value of $m$,
the Generalized CuSum performs substantially worse than the  
 Vector CuSum.

\section{Conclusion}\label{sec:conclusion}



In this paper we revisit the problem of sequential change diagnosis and we propose a novel procedure that enjoys the various  statistical and practical advantages of  existing approaches in the literature.  We analyze  it theoretically and propose a  method for  its design that does not require the estimation of (small) probabilities, nor the simulation of quantities over multiple change-points.  Both this  theoretical analysis and this practical design are formulated  for a more general family of sequential change diagnosis procedures in the literature.

A natural follow-up to this work is an investigation into our conjecture that  the proposed scheme guarantees uniform error control over all possible change-points. Other directions include extensions of the proposed method to non-iid setups and/or composite (non-discrete)  post-change scenarios, as well as  the application of the   proposed adaptive window to the Generalized CuSum algorithm. Finally, the numerical studies  in this work seem to support the use of a pure sequential change detection  algorithm, namely the min-CuSum,  for the sequential change-diagnosis problem, complementing  certain existing works  in the literature (\cite{warnercusum, han2007detection}). This is a topic we plan to consider in more detail in the future. 
%
%


\begin{appendices}

\section{Properties of CuSum Statistics}\label{app:lemmas}

 Throughout all following appendices we use the following notation for the sums of log-likelihood ratios for $n \in \mathbb{N}$: 
 
 \begin{align}\label{def:Z}
 \begin{split}
 Z_i(n) &\equiv \sum_{k=1}^n \ell_i(k), \\
 Z_{ij}(n) &\equiv \sum_{k = 1}^n \log \ell_{ij}(k).
 \end{split}
 \end{align}
 and we further denote partial sums by 
  \begin{align}\label{def:Z_partial}
  \begin{split}
 Z_i(n,m) &\equiv Z_i(n) - Z_i(m), \\
 Z_{ij}(n,m) &\equiv  Z_{ij}(n) - Z_{ij}(m).
 \end{split}
 \end{align}
 
 By defining the above quantities, we may conveniently express the CuSum in the two following ways,

\begin{equation}
 Y_i(n) = \max_{0 \leq k \leq n} Z_i(n,k), 
 \end{equation}
and 
\begin{equation} \label{eqn:cusum_min_form} 
Y_i(n) = Z_i(n) - \min_{0 \leq k \leq n} Z_i(k), 
\end{equation}
with similar expressions obtained for $ Y_{ij}$ by replacing $ Z_i $ with $ Z_{ij}$.

\begin{lemma} \label{lem:vector_cusum_smaller_than_matrix}
 For each $i, j \in [K]$ with $i \neq j$, 
\begin{equation}\label{eqn:vector_cusum_smaller_than_matrix}
 Y_i(n)-Y_j(n)\leq Y_{ij}(n), \quad \forall \; n \in \bN. 
 \end{equation}
\end{lemma}

\begin{IEEEproof}
By \eqref{eqn:cusum_min_form} we have 
\begin{align*}
 Y_i(n)-Y_j(n)
&= \max_{0 \leq k \leq n} Z_i(n,k) - \max_{0 \leq k \leq n} Z_j(n,k) \\
& \leq\max_{0 \leq k \leq n} Z_{ij}(n,k)= Y_{ij}(n).
\end{align*} 
\end{IEEEproof}

\begin{lemma} \label{lem:new}
For each $i,j \in [K]$ with $i \neq j$, for all $x \geq 0$ and  $n \in \mathbb{N}$  we have 
$$ \Pro_{\infty}(Y_i (n) \geq x) \leq e^{-x}.$$
If also $\psi_{ij}$ exists around zero and
\begin{equation}\label{assume_negative_drift}
\Exp_\infty[ \ell_{ij}(1)] <0,
\end{equation} then
$$ \Pro_{\infty}(Y_{ij} (n) \geq x) \leq e^{-r_{ij} x}$$
where  $r_{ij}$ is the positive  root of $ \psi_{ij}.$
\end{lemma}

\begin{IEEEproof} 
For all $x > 0$ and  $n \in \mathbb{N}$ 
\begin{align*}
 \Pro_{\infty}(Y_i (n) \geq x)
&= \Pro_{\infty} \left( \max_{0 \leq m < n} Z_i(n,m) \geq x \right) \\
&= \Pro_{\infty} \left( \max_{0 < m \leq n} Z_i(m) \geq x \right) \leq e^{-x},
\end{align*} 
where the second equality follows by the random walk property and the inequality holds by Ville's supermartingale inequality and the fact that  $\{ \exp [ Z_i(n) ], n \in \bN\}$ is a $\Pro_{\infty}$-martingale with mean 1.

By \eqref{assume_negative_drift} it follows  that 
 $ \left\lbrace \exp [ r_{ij} Z_{ij}(n) ], n \in \mathbb{N}  \right\rbrace$ is a positive martingale with expectation 1 under $\Pro_{\infty}$,
where  $ r_{ij}$  is the positive root of $\psi_{ij}$,  and observe that 
 \begin{align*}
 \Pro_{\infty} \left( Y_{ij}(n) \geq x \right) 
 &= \Pro_{\infty} \left( \max_{0 \leq s < n} Z_{ij}(n,s) \geq x \right) \\
 &= \Pro_{\infty} \left( \max_{0 < s \leq n} Z_{ij}(s) \geq x \right) \leq e^{-r_{ij}x},
\end{align*} 
where the second equality follows from the random walk property the inequality follows from Ville's supermartingale inequality.\\
\end{IEEEproof}

For the next lemma, we introduce  the first  regeneration time of  $Y_{i}$ for each $ i \in [K]$:
\begin{equation} \label{def: eta}
\eta_i \equiv  \inf\{n \in \bN:Y_i(n) = 0\}
\end{equation}
and we note that 
\begin{equation} \label{def: a}
a_i \equiv \sup_{\theta \in (0,1)} |\psi_i (\theta)|>0
\end{equation}
for each $i \in [K]$, where $ \psi_i$ is the cumulant generating function of $ \ell_i(1)$ under $ \Pro_{\infty}$, i.e., 
$$  \psi_{i}(\theta) \equiv \log \left( \Exp_{\infty} \left[ \exp [ \theta \, \ell_{i}(1) ] \right] \right), \quad \theta \in \bR.$$

\begin{lemma} \label{lemma: regeneration bounded by expo}
For each $i \in [K]$, the first regeneration time $\eta_i $ is small under $\Pro_\infty$ in the sense that $$ \Pro_{\infty}(\eta_i > n) \leq e^{-a_in} \quad \forall n \in \mathbb{N},$$ and consequently 
 \begin{equation}\label{eqn:moment_bound}
  \Exp_{\infty}[(\eta_i)^m] \leq \frac{e^{a_i}}{a_i^{m}} m! \quad \text{for all} \; m \in \bN. 
\end{equation} 
 \end{lemma}
 
 \begin{IEEEproof}
For  every $n \in \bN$,   \begin{align*}
\lbrace \eta_i > n \rbrace &= \lbrace Y_i(m) > 0 \;  \text{for all} \; m \in [n]  \rbrace \\
&= \lbrace Z_i(m) > 0 \;  \text{for all} \;m \in [n]  \rbrace \\
&\subseteq \lbrace Z_i(n) > 0 \rbrace. 
\end{align*}
The first equality follows from definitions. The second equality can be easily verified by writing $ Y_i(m) $ in the form \eqref{eqn:cusum_min_form}, and the inequality is obvious. 
Since $ \lbrace \ell_i(n), \, n \in \mathbb{N} \rbrace$ is an iid sequence with negative  expectation  under $\Pro_\infty$ and $ \psi_i(1) = 0$, by  the Chernoff bound we have 
$$ \Pro_{\infty}(\eta_i > n) \leq \Pro_{\infty}(Z_i(n) > 0) \leq e^{-a_i n} \quad \forall \; n \in \bN.$$ 

This condition implies that $\eta_i$ is stochastically smaller than a Geometric random variable with success probability~$ 1 - e^{-a_i}.$ Specifically, if $G$ is distributed as Geom$(1-e^{-a_i})$, then $$\Pro_{\infty}(\eta_i > n) \leq \Pro_{\infty}(G > n), \text{ for all } n \in \bN_0,$$ from which it follows that $$\Exp_{\infty}[\eta_i^m] \leq  \Exp_{\infty}[G^m] \text{ for all } m \in \bN. $$ Moreover, by routine calculation, \begin{align*}
\Exp_{\infty}[G^m] &= m \, \int_0^{\infty} t^{m-1} \Pro_{\infty}(G > t) \, dt \\
 &= m \, \int_0^{\infty} t^{m-1}   e^{-a_i \lfloor t \rfloor} \, dt \\
  &\leq m \, \int_0^{\infty} t^{m-1}   e^{-a_i(t-1)} \, dt \\
  &= \frac{m e^{a_i}}{a_i^m} \int_0^{\infty} u^{m-1} e^{-u} du = \frac{e^{a_i} m!}{a_i^m},
\end{align*} 
and in particular,  $$ \Exp_{\infty}[\eta_i^m] \leq \frac{e^{a_i} m!}{a_i^m}.$$

 \end{IEEEproof}


\section{Properties of CuSum Stopping Times}
Throughout this Appendix, we fix $b, h>0$ and   ${i, j \in [K]}$ such that $i \neq j$.  Moreover, we introduce the following quantity:
$$ \omega_{ij} \equiv \sup_{t \geq 0} \Exp_j \left[ \ell_{ij}(1) - t | \ell_{ij}(1) \geq t \right].$$
That is,  $\omega_{ij}$ is a bound for the expected overshoot of the random walk $Z_{ij}$ above a threshold under $\Pro_j$. 

We recall from Section \ref{sec:CUSUM} that $\sigma_{i}(b)$ is the CuSum stopping time,  with threshold $b$,   for detecting a change from $f$ to $g_i$,
and we introduce the CuSum stopping time,  with threshold $h$,   for detecting a change from $g_j$ to $g_i$, i.e., 
  \begin{align*}
 \sigma_{ij}(h) &\equiv  \inf \lbrace n \in \bN : Y_{ij}(n) \geq h  \rbrace. 
\end{align*}
We denote by  $U_{ij}(x;h)$ the expectation of  $\sigma_{ij}(h)$  under $\Pro_j$ when $ Y_{ij}$ is initialized from  some $x \geq 0 $, i.e., 
\begin{equation}\label{U}
U_{ij}(x;h) \equiv \Exp_j[\sigma_{ij}(h) \,|\, Y_{ij}(0) = x],  \quad x   \geq 0.
\end{equation}

Clearly   $U_{ij}(\cdot; h)$ is a non-increasing function, so that 
$$ U_{ij}(x;h)  \leq U_{ij}(0;h) \equiv \Exp_j[\sigma_{ij}(h)], \quad  \text{ for all } x \geq 0.$$

By  standard properties of the CuSum test
 (see \cite[Section 8.2.6]{tartakovsky2014sequential}) it follows  that 
\begin{align} 
U_{ij}(0;h) =\Exp_j[\sigma_{ij}(h)]  &\geq e^h, \label{ARL_CUSUM_LB} 
\end{align}
and that
\begin{align} \label{u_ineq}
U_{ij}(x;h)  &\geq u_{ij}(x; h),  \quad  x \in [0,h],
\end{align}
where the function $u_{ij}$ is defined as follows:
\begin{equation} \label{def:u}
u_{ij}(x ;h) \equiv  \frac{x - e^{-(h - x)}(h+ \omega_{ij})}{I_{ji}} + (1 - e^{-(h - x)})  \, U_{ij}(0; h) , \quad x\geq 0. 
 \end{equation}
Note that $u_{ij}$ is a differentiable function with derivative 
\begin{equation} \label{derivative2}
 -u_{ij}'(x ;h) = \frac{-1}{I_{ji}} + \frac{h + \omega_{ij}}{I_{ji}}e^{x-h} + U_{ij}(0;h) e^{x-h}. 
 \end{equation}
 Moreover, it follows that \begin{equation}\label{eqn:upper_bound_on_uprime}
 -u_{ij}'(x ;h) \leq \left( \frac{h + \omega_{ij}}{I_{ji}} + U_{ij}(0;h)  \right) \,  e^{x-h}, \quad x \in [0, h]. 
 \end{equation}
For $h\geq 1$, 
\begin{equation} \label{inequality}
 U_{ij}(x ; h) \mathds{1} ( \lbrace x  \leq h \rbrace) \geq u_{ij}(x;h) \quad \text{ for all } x \geq 0. 
 \end{equation}
Indeed, since   $ -e^x \leq -(x + 1)$ for every $x \in \mathbb{R}$, we have 
$$ u_{ij}(x;h) \leq \frac{(x - h)(1 - h)}{I_{ji}}, \quad x \geq h,$$
 and, consequently, for $h \geq 1$, 
 $$ u_{ij}(x;h) \leq 0 \quad \text{ for all } x \geq h.$$

\section{}\label{app:delay}

\begin{IEEEproof}[Proof of Lemma \ref{lem:new_lemma_for delay}]
Fix $\delta>0$. Set
 $$D(b,h) \equiv  \max \lbrace b / I_{i}, h / I^*_i \rbrace(1 + \delta),$$
and observe that 
\begin{align} \label{bound}
\Exp_i &[\tau_i(b,h)] = \sum_{n = 0}^{\infty} \Pro_i \left( \tau_i(b,h) > n \right) \nonumber\\
&\leq  \lceil D (b,h) \rceil + \sum_{n \geq \lceil D (b,h)  \rceil} \Pro_i \left( \tau_i(b,h) > n \right).
\end{align} 

For every $n \in \bN$ we have 
\begin{align*}
& \left\lbrace \tau_i(b,h)> n \right\rbrace \\
&= \bigcap_{m = 1}^n \left\lbrace  Y_i(m) < b \quad \text{or} \quad \min_{j \in [K]: j \neq i} W_{ij}(m) < h \right\rbrace \\
&\subseteq \left\lbrace  Y_i(n) < b \quad \text{ or } \quad \min_{j \in [K]: j \neq i} W_{ij}(n) < h \right\rbrace \\
&= \bigcup_{\stackrel{j \in [K]}{j \neq i}} \left\lbrace W_{ij}(n) < h \right\rbrace \cup \left\lbrace Y_i(n) < b \right\rbrace.
\end{align*}
For $n>D(b,h) $ we have 
$$ b < n \frac{I_{i}}{1 + \delta} \quad \& \quad h < n\frac{I_i^*}{1 + \delta} \leq n\frac{I_{ij}}{1 + \delta} $$
 for every $j \in [K]$ such that $j\neq i$, and consequently, 
\begin{align*}
& \left\lbrace \tau_i(b,h) > n \right\rbrace \\
& \quad \quad \subseteq \bigcup_{\stackrel{j \in [K]}{j \neq i}} \left\lbrace W_{ij}(n)  < n \frac{ I_{ij} }{1 + \delta}  \right\rbrace \cup \left\lbrace Y_i(n) <  n \frac{ I_{i}}{1 + \delta}  \right\rbrace
\end{align*}
 for every  $n > D(b,h)$. Combining this with \eqref{bound} we obtain
 \begin{align} \label{bound2}
 \begin{split}
 \Exp_{i}[\tau_i(b,h)] &\leq  \lceil D(b,h) \rceil \\
&+ \sum_{n =1 }^\infty \Pro_i \left( Y_i(n) < n I_{i} (1 + \delta)^{-1} \right)  \\
& + \sum_{n =1 }^\infty \sum_{\stackrel{j \in [K]}{j \neq i}} \Pro_i \left( W_{ij}(n) < n I_{ij} (1 + \delta)^{-1}  \right).
\end{split}
\end{align} 
The third term in the upper bound converges by assumption. Moreover,  since  $Y_i(n) \geq Z_i(n)$ for every 
$ n \in \bN$ and   by the Chernoff bound it follows that $\Pro_i \left( Z_{i}(n) \leq\rho n  \right)$ is an exponentially decaying sequence for every  $\rho < I_i$, we conclude that 
 $\Pro_i \left( Y_{i}(n) \leq \rho n  \right)$ is also an exponentially decaying sequence for every  $\rho < I_i$. Therefore,  the second series in the upper bound converges, and this completes the proof. \\
 \end{IEEEproof}


\begin{IEEEproof}[Proof of Lemma \ref{lem:det_delay_bound}]
For every $n \in \bN$ we have 
\begin{align*} 
\min_{R_i(n) \leq k \leq n} Z_{ij}(k)
&\leq  Z_{ij}(R_i(n))\\
&=\sum_{m = 1}^{R_i(n) }\ell_{ij}(m)\leq\sum_{m = 1}^{R_i(\infty) } |\ell_{ij}(m)|,
\end{align*}
where $R_i(\infty) $ is the final regeneration time of $Y_i$, i.e., 
  $$R_i(\infty) \equiv \sup \lbrace n \in \bN
  : Y_{i}(n) = 0 \rbrace, $$
and $Z_i$ is defined in \eqref{def:Z}. Consequently,
\begin{align*} 
Y^{\prime}_{ij}(n) &= \max_{R_i(n) \leq k \leq n} Z_{ij}(n,k) \\
&= Z_{ij}(n) - \min_{R_i(n)  \leq k \leq n} Z_{ij}(k)\\
&\geq Z_{ij}(n)- \sum_{m = 1}^{R_i(\infty) } |\ell_{ij}(m)|.
\end{align*}

  Let $\rho \in (0, I_{ij})$. Then  there is a $\delta>0$ such that $\rho+\delta< I_{ij}$, and as a result    we have 
\begin{align*}
\Pro_i \left( Y^{\prime}_{ij}(n) \leq \rho  n  \right) &\leq  \Pro_i \left( Z_{ij}(n) \leq (\rho + \delta ) n  \right)\\
&+   \Pro_i \left( \sum_{m=1}^{R_i(\infty) }| \ell_{ij}(m)| \geq n \delta   \right).
\end{align*}
By the Chernoff bound, the first term in the upper bound goes to 0 exponentially fast in $n$.  Thus, it remains to show that the second term is summable, for which it suffices to show that 
$$\Exp_i \left[ \sum_{m=1}^{R_i(\infty) } |\ell_{ij}(m)|  \right] < \infty.$$
By the first assumption of the lemma we have 
$$ || \ell_{ij}(1) ||_p \equiv \sqrt[p]{\Exp_i \left[ |\ell_{ij}(1)|^p \right]} <\infty,$$ 
and as a result
\begin{align*}
\Exp_i \left[ \sum_{m=1}^{R_i(\infty) } |\ell_{ij}(m)|  \right] &=     \sum_{m = 1}^{\infty} \mathbb{E}_i \left[ |\ell_{ij}(m)| \;  \mathds{1} ( \lbrace R_i(\infty)  \geq m \rbrace ) \right] \\
&\leq    || \ell_{ij}(1) ||_p  \;    \sum_{m= 1}^{\infty} 
\sqrt[q]{\Pro_i \left( R_i(\infty)  \geq m \right) }
 \end{align*}
 where $ q \equiv  p / (p-1), $ the equality follows by Tonelli's theorem, and  the inequality by  H{\"o}lder's inequality. Therefore, it suffices to show that the series in the upper bound converges.  By the  second assumption of the lemma there is an  $\epsilon>0$ such that 
$$ \Exp_i \left[ \left( \ell_{i}^-(1) \right)^{2 + \epsilon +q} \right] < \infty,$$
and by \cite[Theorem 1]{janson1986moments}
it follows that
$$\Exp_i[(R_i(\infty) )^{1 + \epsilon +q}] < \infty,$$
or equivalently
$$\sum_{m= 1}^{\infty} m^{\epsilon +q} \, \Pro_i(R_i(\infty)  \geq m) < \infty.$$  As a result,  there is some  $ s \in \bN $ such that  
 $$ m^{\epsilon + q} \, \Pro_i(R_i(\infty)  \geq m) < 1 \quad \text{for all} \quad   m \geq s,$$
and  this further implies  that 
 $$ \sum_{m = s}^{\infty} \sqrt[q]{\Pro_i \left(  R_i(\infty)  \geq m \right) } \leq   \sum_{m = s}^{\infty}  m^{-(1 + \epsilon /q)} <\infty,$$
 which completes the proof. \\
\end{IEEEproof}

\section{}\label{app:false_alarm_bound}

\begin{IEEEproof}[Proof of Lemma \ref{lem:false_alarm_bound}]
We fix $b, h>0$  and  observe that the assumption of the Lemma implies that, for all $i, j \in [K], j \neq i$,
 $$ \sup_{n \in \bN_0} \Pro_{\infty}(W_{ij}(n) \geq h) \leq Q_{ij} \, e^{-q_{ij}h}.  $$
As a result,  for every $i \in [K]$  we have 
 $$ \sup_{n \in \bN_0} \Pro_{\infty}  \left( \min_{j \in [K] : j \neq i} W_{ij}(n) \geq h \right) \leq G_i(h), $$ 
where  
$$ G_i(h) \equiv \min_{j \in [K]: j \neq i}  \left\lbrace Q_{ij} e^{-q_{ij} h} \right\rbrace. $$
Since also 
\begin{align*}
\tau(b,h) &\geq \inf \left\{ n \in \bN : \max_{i \in [K]} \min_{j \in [K]: j \neq i} W_{ij}(n) \geq h \right\},
\end{align*}
for every $n \in \bN$ we have 
\begin{align*}
\Pro_{\infty}(\tau(b,h) \leq n) &\leq \sum_{m = 1}^n \sum_{i = 1}^K \Pro_{\infty} \left( \min_{j \in [K]: j \neq i} W_{ij}(m) \geq h \right) \\
&\leq \sum_{m = 1}^n \sum_{i = 1}^K  G_i(h) \equiv  n G(h),
\end{align*} 
and, consequently,
 $$ \Pro_{\infty}(\tau(b,h) > n) \geq \left( 1 - n  G(h)\right)^+. $$
As a result, 
\begin{align*}
\Exp_{\infty} \left[ \tau(b,h) \right] &= \sum_{n = 0}^{\infty} \Pro_{\infty}(\tau(b,h) > n) \\
&\geq  \sum_{n = 0}^{\infty}\left( 1 - n  G(h) \right)^+\\
&=\sum_{n = 0}^{\lfloor 1/G(h) \rfloor} \left( 1 - n  G(h) \right)  \geq \frac{1}{2G(h)},
\end{align*}
which is what we wanted to show. \\
\end{IEEEproof}

\section{}\label{app:uniform_exponential_bound}

\begin{IEEEproof}[Proof of Theorem \ref{thm:uniform_exponential_bound}]
Throughout the proof, we fix arbitrary  $x>0 $ and  $n \in \bN_0  $. Since, by definition,  $Y^{\prime}_{ij} \leq Y_{ij}$, we have 
 \begin{align*}
 \Pro_{\infty} \left( Y^{\prime}_{ij}(n) \geq x \right)  \leq \Pro_{\infty} \left( Y_{ij}(n) \geq x \right), 
\end{align*} 
 and the result follows from Lemma \ref{lem:new} when   $\Exp_\infty[ \ell_{ij}(1)]<0$.  Thus, in the remainder of the proof we focus on the case that  $\Exp_\infty[ \ell_{ij}(1)] \geq 0$. By the definition of $Y'_{ij}$ and the total probability law, for any $ \delta > 0 $ we have  
 \begin{align}  \label{very_first_bound}
 \begin{split}
\Pro_{\infty} \left( Y^{\prime}_{ij}(\nu) \geq x \right)  & \leq \Pro_{\infty} \left( \max_{\nu - \delta x < s \leq \nu} Z_{ij}(\nu, s)   \geq x \right)  \\
& \quad + \Pro_{\infty}\left( \nu - R_i(\nu) \geq \delta x \right). 
\end{split}
\end{align}
We start by upper bounding the first term in the right-hand side of  \eqref{very_first_bound}. Since ${ \{ Z_{ij}(n), n \in \bN \} }$ is a random walk under $\Pro_\infty$, 
\begin{align*}  
\Pro_{\infty} \left(  \max_{\nu - \delta x < s \leq \nu} Z_{ij}(\nu,s)   \geq x \right)
&= \Pro_{\infty} \left( \max_{1 \leq s \leq \delta x} Z_{ij}(s) \geq x \right) .
\end{align*}
By assumption,  $ \psi_{ij} $ is finite around zero, so  $\psi_{ij}(\theta)<\infty$ for $\theta>0$ small enough. For any such $\theta$,
 $$\{\exp [ \theta Z_{ij}(n) ], n \in \bN\}$$
  is a positive submartingale under $\Pro_{\infty}$, as can be seen by combining Jensen's inequality with the fact that $ \Exp_{\infty}[\ell_{ij}(1)] \geq 0$, and by Doob's submartingale inequality it follows that
\begin{align}  \label{first_bound_old}
\Pro_{\infty} \left( \max_{1 \leq s \leq \delta x} Z_{ij}(s) \geq x \right)   \leq \exp [- \theta x +  \lfloor \delta x \rfloor  \psi_{ij}(\theta)].
\end{align}

Since $ {\Exp_{\infty}[\ell_{ij}(1)] \geq 0}$, by the strict convexity of $\psi_{ij}$ it follows that 
$\psi_{ij}(\theta) > 0$ for  $ \theta > 0$ small enough. Therefore, since also $ \lfloor x \delta \rfloor > 0$, we can rewrite \eqref{first_bound_old} as
\begin{equation}\label{first_bound}
 \Pro_{\infty} \left( \max_{1 \leq s \leq \delta x} Z_{ij}(s) \geq x \right) \leq \exp [- ( \theta  -   \delta \psi_{ij}(\theta)) \, x].
 \end{equation}

We continue with the second term in the right-hand side of  \eqref{very_first_bound}.  By Markov's inequality, for any $\theta > 0$, 
\begin{align*}
& \Pro_{\infty}\left( \nu - R_i(\nu) \geq \delta x \right) \\
&\leq e^{-\theta \delta x} \, \Exp_{\infty}  \left[ e^{\theta ( \nu - R_i(\nu))} \right] \\
&= e^{-\theta \delta x} \,  \sum_{m = 0}^{\infty} \frac{ \theta^m}{m!} \;  \Exp_{\infty}[ (\nu - R_i(\nu))^{m}],
\end{align*}
where the equality follows from a Taylor series expansion and the Monotone Convergence Theorem.

Since $ \nu - R_i(\nu)$ is the age at time $\nu$ of a discrete renewal process whose events are the regenerations of $ \lbrace Y_i(n), n \in \mathbb{N} \rbrace$, by Lorden's inequality (see \cite[Theorem 3]{lorden1970excess} and the comments therein), 
 for every $m \in \bN$ and $\nu \in \bN$  we have:
 $$\Exp_{\infty} \left[ (\nu - R_i(\nu))^m \right] \leq \frac{m+2}{m+1} \; \frac{\Exp_{\infty} [\eta_i^{m+1}]}{\Exp_{\infty}[\eta_i]}, 
 $$
 where $\eta_i$ is the first  regeneration time of the sequence ${ \{Y_{i}(n), n \in \bN \} }$, defined in \eqref{def: eta}.
By Lemma \ref{lemma: regeneration bounded by expo}  we then conclude that 
$$ \Pro_{\infty}(\nu - R_i(\nu) \geq \delta  x) \leq \frac{e^{-\theta \delta x} e^{a_i}}{ \Exp_{\infty} [\eta_i] a_i} \, \sum_{m=0}^{\infty} \left(\theta/ a_i \right)^m (m+2), $$
where $a_i$ is defined in \eqref{def: a}, which means that, for every $ \theta < a_i$, there is a constant $K_{\theta}>0$, which does not depend on $x$ or $\delta$, such that 
\begin{equation} \label{second_bound}
\Pro_{\infty}(\nu - R_i(\nu) \geq \delta x) \leq K_{\theta} \, e^{-\theta \delta  x}.
\end{equation}

Combining \eqref{first_bound} and \eqref{second_bound}  we conclude that, for every ${\theta \in (0, a_i)}$ such that $\psi_{ij}(\theta)<\infty$,
$$\Pro_{\infty} \left( Y^{\prime}_{ij}(\nu) \geq x \right)  \leq e^{-x(\theta - \delta \psi_{ij}(\theta))} + K_{\theta}\,  e^{-\theta \delta x}.$$
The above inequality implies the desired result once we select a sufficiently small $\theta \in (0, a_i)$  so that 
$\psi_{ij}(\theta)<\infty$ and then a sufficiently small $\delta>0$ so that  $\theta - \delta \psi_{ij}(\theta) > 0$, or equivalently  $ \delta \in (0, \theta / \psi_{ij}(\theta) )$, recalling that $ \psi_{ij}(\theta) > 0 $ for $ \theta >0$ by strict convexity. \\
\end{IEEEproof}

\section{}\label{app:family_misspecification_bound_theorem}

We present a lemma to help us in the proof of the Theorem. 
\begin{lemma}\label{lem:prob_to_expectation}
Let $ T_1, T_2 $ be $ \bN$-valued random variables on a probability space with a probability measure $ \Pro $ and finite expectation $ \Exp[T_i] < \infty, i = 1,2.$ If $ \Pro(T_1 > T_2 ) > 0, $ then $$ \Pro \left( T_1 > T_2\right) = \frac{\Exp \left[ T_1 - (T_2 \wedge T_1) \right] }{\Exp \left[ T_1 - T_2  \, | \, T_1 > T_2\right]}. $$
\end{lemma}
\begin{IEEEproof}
It suffices to observe that
\begin{equation}
\Exp \left[ T_1 - (T_2 \wedge T_1) \right] = \Exp \left[ T_1 - (T_2 \wedge T_1)\,  | \, T_1 > T_2 \right] \times \Pro(T_1 > T_2)
\end{equation} and also that $$
\Exp \left[ T_1 - (T_2 \wedge T_1)\,  | \, T_1 > T_2) \right] = \Exp \left[ T_1 - T_2 \,  | \, T_1 > T_2 \right].
$$
\end{IEEEproof}

Now, proceeding with the proof of the Theorem, we follow similar steps as in the proof of \cite[Theorem 2]{nikiforov2000simple}.

\begin{IEEEproof}[Proof of Theorem \ref{th:family_misspecification_bound_theorem}]
Throughout this proof, for simplicity we write $ \tau_i, \tau, \widehat{\tau}$,  instead of  $\tau_i(b,h),  \tau(b,h)$, $\widehat{\tau}(b,h)$, respectively, but it is important to keep this dependence in mind. By the definition of $\tau$ and  $\widehat{\tau}$ in \eqref{family}, 
$$\{\widehat{\tau}=i\}  \subseteq 
\{\tau_{i} \leq \tau_j \}.
$$
Moreover, since  $\tau\leq \tau_i$,  on the event $\{\tau>\nu\}$ we have 
\begin{align*}
 \tau_i  &= \inf \lbrace n > \nu : Y_i \geq b, \,  W_{ik}(n) \geq  h, \; \;   \forall \,  k \in [K], k \neq i \rbrace \\
 &\geq \inf \lbrace n > \nu : W_{ij}(n) \geq h \rbrace\\
  &\geq \inf \lbrace n > \nu : W^\nu_{ij}(n) \geq h \rbrace = \tau^\nu_{ij} ,
 \end{align*} 
 where the last inequality holds because, by
the assumption that  \eqref{pathwise bound} holds  and the definition of $W_{ij}^\nu$  in \eqref{W_nu}, 
\begin{align} \label{implication}
 W_{ij} (n) \leq W^{\nu}_{ij}(n) \quad \text{for all} \quad n \in \bN,
 \end{align}
 and the last equality is simply the definition of  $\tau^\nu_{ij}$ in \eqref{tau_nu}.  Therefore,
\begin{align*} 
\{\widehat{\tau}=i, \tau>\nu \} 
&\subseteq 
\{ \tau_i \leq \tau_j, \tau>\nu \} \subseteq  \{ \tau^\nu_{ij}  \leq \tau_j, \tau>\nu \},
\end{align*}
and consequently
\begin{align}\label{eqn:1111}
\Pro_{\nu, j}( \widehat{\tau}=i | \; \tau > \nu)  
&\leq \Pro_{\nu, j} ( \tau_{ij}^{\nu} \leq  \tau_j |\;  \tau > \nu).
\end{align}
Using  assumption  \eqref{condition technical} of Theorem \ref{th:family_misspecification_bound_theorem}  and applying
Lemma \ref{lem:prob_to_expectation} with  $\Pro\equiv  \Pro_{\nu, j}(\cdot \, | \,  \tau > \nu)$, 
we  have
\begin{align}\label{eqn:1}
\begin{split}
 & \Pro_{\nu, j} ( \tau_{ij}^{\nu} \leq  \tau_j  \, |\,  \tau > \nu) 
= 1- \Pro_{\nu, j} ( \tau_{ij}^{\nu} > \tau_j  \,| \, \tau > \nu)\\
&=1-   \frac{\Exp_{\nu, j}[\tau^{\nu}_{ij} - \tau_j \wedge \tau^\nu_{ij}  \, | \, \tau > \nu]}{\Exp_{\nu, j}[\tau^{\nu}_{ij} -  \tau_j \,  |\,  \tau > \nu, \tau^{\nu}_{ij} >  \tau_j]}.
\end{split}
\end{align}
We proceed by further lower bounding the denominator.
Since $ \tau \leq \tau_j, $ on the event $\lbrace \tau > \nu, \tau^{\nu}_{ij} > \tau_j \rbrace$ we have 
\begin{itemize}
\item $\nu < \tau \leq \tau_j < \tau_{ij}^{\nu},$
\item $\tau_{ij}^{\nu} - \tau_j $ is a function of $W_{ij}^{\nu}(\tau_j)$ and $ 
\{X_{n}, n > \tau_j \}$
\end{itemize}
and, as a result, 
\begin{equation}
\Exp_{\nu, j}[ \tau_{ij}^{\nu} - \tau_j  \, |\, \tau > \nu, \tau_{ij}^{\nu} > \nu, \mathcal{F}_{\tau_j}] = U_{ij}(W_{ij}^{\nu}(\tau_j) ; h) \cdot \mathds{1} ( \{\tau > \nu, \tau_{ij}^{\nu} > \nu\} ) 
\end{equation} 
where $U_{ij}$ is defined as in \eqref{U}.
Since  $ W_{ij}^{\nu}(\tau_j) < h $ on $ \lbrace \tau_{ij}^{\nu} > \tau_j > \nu \rbrace$ and  the function  $U_{ij}(\cdot\,  ; h) $ is  non-increasing, we further obtain 
\begin{align*}
\Exp_{\nu, j}[ \tau_{ij}^{\nu} - \tau_j  \, |\, \tau > \nu, \tau_{ij}^{\nu} > \nu, \mathcal{F}_{\tau_j}] 
&\leq U_{ij}(0; h)= \Exp_j[\sigma_{ij}],
\end{align*} 
and,  by the law of iterated expectation,  we conclude that:
\begin{equation}\label{eqn:3}
\Exp_{\nu, j}[ \tau_{ij}^{\nu} - \tau_j  \, | \, \tau > \nu, \tau_{ij}^{\nu} > \nu] \leq \Exp_j[\sigma_{ij}].
\end{equation}

Combining \eqref{eqn:1111}-\eqref{eqn:3}, we obtain
\begin{align*}
\Pro_{\nu, j}&(\widehat{\tau} = i | \tau > \nu) \leq 1 - \frac{\Exp_{\nu, j}[\tau^{\nu}_{ij} - \tau_j \wedge \tau^\nu_{ij}  | \; \tau > \nu]}{\Exp_{j}[\sigma_{ij}]} \\
&\leq \frac{\Exp_{j}[\sigma_{ij}]  -  \Exp_{\nu, j}[\tau^{\nu}_{ij} - \nu  | \; \tau > \nu]}{\Exp_{j}[\sigma_{ij}]} +  \frac{\Exp_{\nu, j}[ \tau_j - \nu | \; \tau > \nu]}{\Exp_{j}[\sigma_{ij} ]},
\end{align*}
and applying  \eqref{ARL_CUSUM_LB} we arrive at
\begin{align*}
\Pro_{\nu, j}(\widehat{\tau} = i | \tau > \nu) 
&\leq  \frac{\Exp_{j}[\sigma_{ij}]  -  \Exp_{\nu, j}[\tau^{\nu}_{ij} - \nu  | \; \tau > \nu]}{\Exp_{j}[\sigma_{ij}]} \\
&+   e^{-h} \, \Exp_{\nu, j}[ \tau_j  - \nu | \; \tau > \nu]. 
\end{align*}
Comparing with \eqref{show}, it is clear that it  suffices to show that 
$$ \frac{\Exp_{j}[\sigma_{ij}]  -  \Exp_{\nu, j}[\tau^{\nu}_{ij} - \nu  | \; \tau > \nu]}{\Exp_{j}[\sigma_{ij}]} \leq  \frac{C_{ij}}{1 - c_{ij}}e^{-c_{ij}h} \left( 1 + \phi_{ij}(h) \right),
$$
  where $\phi_{ij}$ is a function that goes to $0$ as $h \to \infty$. 
In order to do so, we add and subtract  $u_{ij}(0;h)$ (defined in \eqref{def:u}) on the left-hand side, 
 \begin{align}\label{eqn:star1}
 \begin{split}
  \Exp_j[\sigma_{ij}] -  &\Exp_{\nu, j}[\tau^{\nu}_{ij} - \nu| \tau > \nu]  \\
 = &\;  U_{ij}(0;h) - u_{ij}(0;h)  +   u_{ij}(0;h) - \Exp_{\nu, j}[\tau^{\nu}_{ij} - \nu\,| \, \tau > \nu].
 \end{split} 
 \end{align} 
 We also obtain from \eqref{def:u} that 
 \begin{align} \label{eqn:star2}
 U_{ij}(0;h) - u_{ij}(0;h) = e^{-h} \left( \frac{h + \omega_{ij}}{I_{ji}} + U_{ij}(0;h) \right). 
 \end{align}
Moreover, since $ \tau^{\nu}_{ij} - \nu$  is a function of only $W_{ij}(\nu) $ and $ \lbrace X_n, n > \nu \rbrace, $ 
 \begin{align*}
\Exp_{\nu, j} \left[ \tau^{\nu}_{ij} -  \nu  | \cF_{\nu} , \tau>\nu \right] = U_{ij}(W_{ij}(\nu)  ; h) \cdot \mathds{1} ( \lbrace \tau > \nu \rbrace ) 
\end{align*}
and, by the law of iterated expectation, 
\begin{align*}
\Exp_{\nu, j}[\tau^{\nu}_{ij} - \nu \, | \,  \tau > \nu ] 
& =  \Exp_{\nu, j}[ U_{ij}(W_{ij}(\nu)  ; h) \, |\,  \tau > \nu].
\end{align*}
We conclude that
\begin{align*}
&  \Exp_{\nu, j}[ U_{ij}(W_{ij}(\nu)  ; h) \, |\,  \tau > \nu] \\
&\geq  \Exp_{\nu, j}[U_{ij}(W_{ij}(\nu) \wedge h ; h) \cdot \mathds{1} ( \lbrace W_{ij}(\nu) \leq h  \rbrace )| \tau > \nu] \\
&\geq \Exp_{\infty}[u_{ij}(W_{ij}(\nu) \wedge h ; h) | \tau > \nu],
\end{align*}
where  the first inequality holds because 
$U(\cdot \, ; h)$ is non-negative   and the last one, when $ h \geq 1$, by 
\eqref{inequality}. Putting these together,  we obtain
\begin{align*}
\begin{split}
&  u_{ij}(0;h)  - \Exp_{\nu, j}[\tau^{\nu}_{ij} - \nu| \tau > \nu] \\
  &\leq \Exp_{\infty}[u_{ij}(0;h) -  u_{ij}(W_{ij}(\nu) \wedge h  ; h) \, | \,  \tau > \nu].
  \end{split}
\end{align*} 
Since the function  $u_{ij}(\cdot \, ;h)$ is differentiable in $[0,h]$,   by the Fundamental Theorem of Calculus and inequality  \eqref{eqn:upper_bound_on_uprime}  we  obtain 
 \begin{align*}
  u_{ij}(0;h) -  u_{ij}(W_{ij}(\nu) \wedge h ; h) &= - \int_0^{h} u'_{ij}(x ; h) \cdot  \mathds{1} ( \lbrace W_{ij}(\nu) \geq x \rbrace )  \; dx   \\
&\leq    \left( \frac{h + \omega_{ij}}{I_{ji}} + U_{ij}(0;h)  \right)  e^{-h} \; \int_0^{h} e^{x} \;  \mathds{1} ( \lbrace W_{ij}(\nu) \geq x \rbrace ) \; dx.
\end{align*}
By Tonelli's theorem  and  assumption \eqref{condition on the change-point} we further obtain
 \begin{align*}
 \Exp_{\infty}\left[\int_0^{h} e^{x} \;  \mathds{1}( \lbrace W_{ij}(\nu) \geq x \rbrace )  \; dx \, | \,  \tau > \nu \right] &\leq \int_0^{h} e^{x} \,  \Pro_{\infty}(W_{ij}(\nu) \geq x \, \big| \, \tau > \nu) \, dx \\
&\leq   \int_0^{h} \ e^{x} \; C_{ij} e^{-c_{ij}x} \; dx \leq  \frac{C_{ij}}{1 - c_{ij}} \, e^{(1 - c_{ij})h}
 \end{align*} 
and, consequently, 
\begin{equation}\label{eqn:star3}
u_{ij}(0;h)  - \Exp_{\nu, j}[\tau^{\nu}_{ij} - \nu| \tau > \nu] \leq \frac{C_{ij}}{1 - c_{ij}} \, e^{-c_{ij}h} \, \left( U_{ij}(0;h) +  \frac{h + \omega_{ij}}{I_{ji}} \right).
\end{equation}
By \eqref{eqn:star1}-\eqref{eqn:star3} we obtain 
 \begin{align}\label{eqn:star4}
 \begin{split}
 &\Exp_j[\sigma_{ij}] -  \Exp_{\nu, j}[\tau^{\nu}_{ij} - \nu| \tau > \nu]  \\
  & \quad \leq   e^{-h} \left( \frac{h + \omega_{ij}}{I_{ji}} + U_{ij}(0;h) \right)  +  \frac{C_{ij}}{1 - c_{ij}} \, e^{-c_{ij}h} \, \left( U_{ij}(0;h) +  \frac{h + \omega_{ij}}{I_{ji}} \right).
  \end{split}
 \end{align} 
 Dividing both sides by $ \Exp_{j}[\sigma_{ij}]$, which is equal by definition to $ U_{ij}(0;h)$, and then applying inequality \eqref{ARL_CUSUM_LB} completes the proof.  \\
\end{IEEEproof}

\begin{IEEEproof}[Proof of Corollary \ref{coro: h_choice}]
By Theorem \ref{thm:delay_bound} it is clear that, for any given  $\alpha \in (0,1)$, the upper bound in \eqref{show} goes to 0 as $h \to \infty$. Therefore, for any  given
$\alpha, \beta \in (0,1)$, there is an $h_{\alpha, \beta}$ so that the worst-case probability of false isolation does not exceed $\beta$. Therefore, to prove the corollary it is enough to show that, for all $ \nu \in N_{\bs{C}}$, 
$$  \Pro_{\nu, j}(\tau_{ij}^{\nu}(h) > \tau_j(b,h) \,|\,  \tau(b,h) > \nu) > 0,$$
at least for $h$ large enough and $ b = b_{\alpha}$, 
where $ \tau_{ij}^{\nu}$ is defined in \eqref{tau_nu}.  In the rest of the proof, we fix  $\nu \in N_{\bs{C}}$,  and for simplicity we suppress $b$ and $h$ in the notation, except where necessary.

First of all, we observe that 
\begin{align} \label{dec}
\begin{split}
\Pro_{\nu, j} ( \tau_{ij}^{\nu} >  \tau_j \,| \, \tau > \nu) &= \Pro_{\nu, j} ( \tau_{ij}^{\nu} >  \tau_j  \, | \,  Y'_{ij}(\nu) \leq  h/2, \tau > \nu )  \\
&\cdot \Pro_{\nu, j} ( Y'_{ij}(\nu) \leq  h/2 \, | \,  \tau > \nu).
\end{split}
\end{align}
Since  $\nu \in N_{\bs{C}}$, by  \eqref{condition on the change-point}  it follows  that 
\begin{equation*}
\Pro_{\infty} (Y'_{ij}(\nu) \leq  h/2 | \tau > \nu) \geq 1  -C_{ij}e^{-c_{ij}h/2},
\end{equation*}
where the lower bound is positive for  $h > 2 \log (C_{ij}) / c_{ij}$.  Therefore, it remains to show that the first factor in \eqref{dec} is positive at least for large $h$.  To this end, we first show that it is bounded below by the probability that 
 $\sigma_{ij}$ is larger than $\tau_j$ when its statistic $Y_{ij}$ is initialized from $h/2$, i.e., 
\begin{align} \label{dec2} 
\Pro_{\nu, j} ( \tau_{ij}^{\nu} >  \tau_j  \, | \,  Y'_{ij}(\nu) \leq  h/2, \tau > \nu )  &\geq   \Pro_{j} \big( \sigma_{ij} >  \tau_j  \, | \,  Y_{ij}(0) = h/2 \big).\end{align}
Indeed, by the law of iterated expectation,
\begin{align}\label{eqn:total_exp}
\begin{split}
&\Pro_{\nu, j} ( \tau_{ij}^{\nu} >  \tau_j  \, | \,  Y'_{ij}(\nu) \leq h / 2, \tau > \nu )  = \\
&\Exp_{\nu,j} \Big[ \zeta(Y'_{ij}(\nu), Y_j(\nu), Y'_{j1}(\nu), \ldots, Y'_{jK}(\nu) \big) \Big|  \tau > \nu, Y'_{ij}(\nu) \leq h / 2 \Big],
\end{split}
\end{align}
where 
\begin{align*}
&\zeta(y'_{ij}, y_j, y'_{j1}, \ldots, y'_{jK}) \\
&\equiv
\Pro_{\nu, j} \big( \tau_{ij}^{\nu} >  \tau_j  \, | \,  Y'_{ij}(\nu)=y'_{ij}, Y_j(\nu)= y_j, Y'_{jk}(\nu)=y'_{jk}, 1 \leq k \leq K \big).
\end{align*}

The stopping time $\tau_{ij}^{\nu}$ depends on~$\cF_\nu$ only through $Y'_{ij}(\nu)$, whereas $\tau_j$ depends on $\cF_\nu$ only through 
$$Y_j(\nu), Y'_{j1}(\nu), \ldots, Y'_{jK}(\nu),$$ 
and $\zeta$ is  increasing in its first argument and decreasing in each of the other arguments. Moreover,  $\tau^\nu_{ij}-\nu$ has the same  distribution as  $\sigma_{ij}$ when the latter is initialized from  $Y'_{ij}(\nu)$. As a result, we conclude that, for all values of $(y'_{ij}, y_j, y'_{j1}, \ldots, y'_{jK})$ in $\lbrace \tau > \nu, Y'_{ij}(\nu) \leq h / 2  \rbrace$, 
\begin{align*}
\zeta(y'_{ij}, y_j, y'_{j1}, \ldots, y'_{jK}) &\geq 
\zeta(h/2,0, \ldots, 0) =   \Pro_{j} \big( \sigma_{ij} >  \tau_j  \, | \,  Y_{ij}(0) = h/2 \big),
\end{align*}
which proves \eqref{dec2}. 


It remains to show that lower bound in \eqref{dec2} is positive. For this,  it suffices to show  that 
\begin{equation}\label{dec3}
\Pro_{j} \big( \sigma_{ij}(h) >  \tau_j(b,h)  \, | \,  Y'_{ij}(0) = h/2 \big) \geq \Pro_j \big( \sigma_{ij}(h/2) >  \tau_j(b,h) \big) 
\end{equation}
since the latter is greater than zero.
Indeed, by  \eqref{ARL_CUSUM_LB} and Theorem \ref{thm:delay_bound} it follows that  $\Exp_j[\sigma_{ij}(h/2)]$ has an exponential lower bound  in~$h$, and $ \Exp_j[\tau_j(b_{\alpha}, h)]$ a linear upper bound in~$b_{\alpha}$ and~$h$  under~$\Pro_j$, which implies that the  lower bound  is positive, at least for~$h$ large enough.  

To prove \eqref{dec3}  it suffices   to show that   $ \sigma_{ij}(h)$ initialized from $Y_{ij}(0) = h/2$ is stochastically  larger than $ \sigma_{ij}(h/2) $  initialized from  $ Y_{ij}(0) = 0.$ Indeed,  the first one corresponds to the  first time $Y_{ij}$ crosses the threshold~$h$ when initialized from $h/2$, whereas  the second  to the first time $Y_{ij}$ crosses the threshold~$h/2$ when initialized from $0$.  In both cases,  $Y_{ij}$ needs to increase by $h/2$,  but in the first case it can fall below its initialization point, whereas in the second it is reflected at its initialization point and, as a result, it is pathwise closer to the threshold than in the first case. 
\end{IEEEproof}

\begin{IEEEproof}[Proof of Theorem \ref{th:optimality}]
By Theorem \ref{thm:delay_bound} it follows that for every $\delta>0$  there is constant $C_\delta>0$ so that 
\begin{align} \label{delay under 00}
\cJ_i \left[ \tau_i(b,h) \right] \leq  (1+\delta)   \left( \max \lbrace b / I_{i}, h / I^*_i \rbrace + C_\delta \right). 
\end{align}

By  this  and  Theorem \ref{th:family_misspecification_bound_theorem}, there are constants $ {C > 0}, {c \in (0,1)}, $ and a function $ \phi(h)$ that goes to zero as $h $ goes to infinity  so that 
\begin{align}\label{eqn:exact_upper_bound_thm8}
\begin{split}
&  \max_{j \in [K]} \sup_{\nu \in N_{\bs{C}}} \Pro_{\nu, j}(\widehat{\tau}(b,h) \neq j | \tau(b,h) > \nu) \\
 &\leq  C \left(   e^{-c h} (1 + \phi(h) )  +   e^{-h}  (\max\{ b_{\alpha}, h \} +  1)  \right) 
\end{split}
\end{align}
Therefore, if we set 
$$
h_\beta \equiv \left( |\log \beta|+\log C \right)/c, \quad \beta \in (0,1),
$$
then by  \eqref{eqn:exact_upper_bound_thm8} it follows that 
$$  \max_{j \in [K]} \sup_{\nu \in N_{\bs{C}}} \Pro_{\nu, j}(\widehat{\tau}(b_\alpha, h_\beta) \neq j | \tau(b,h) > \nu) \lesssim \beta
$$
as $\alpha, \beta \to 0$ so that $| \log \beta |  \gg \log |\log \alpha|$
and by  \eqref{delay under 00} we conclude that 
$$ \mathcal{J}_i[\tau(b_{\alpha}, h_{\beta})] \lesssim 
\frac{|\log \alpha|}{I_i}$$
as $\alpha, \beta \to 0$ so that $ \log |\log \beta| \ll 
\log |\log \alpha|$, which completes the proof. 
\end{IEEEproof}

\newpage
\section{}

Illustration of computing threshold regions: 
\begin{figure}[h!]
\begin{center}	
	
		\subfloat[Proposed]{
			\includegraphics[scale=1]{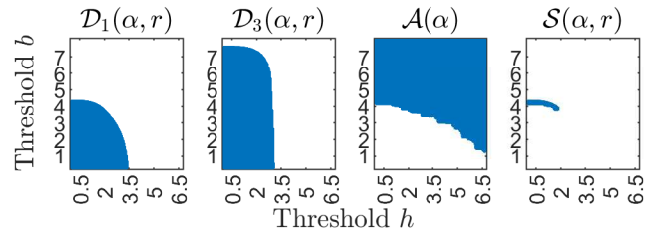}  
			\label{fig:proposed_threshold_regions_1}
		}
		
		\subfloat[Matrix CuSum]{
			\includegraphics[scale=1]{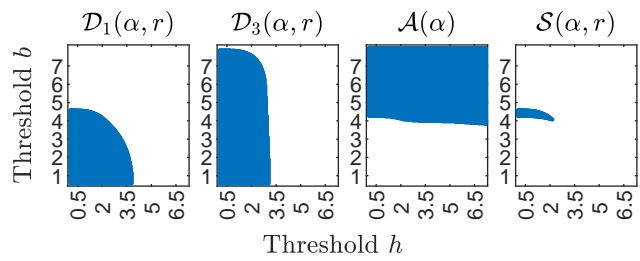}  
			\label{fig:matrix_cusum_threshold_regions_1}
		}
		
		\captionsetup{justification=raggedright, singlelinecheck=false}
		\caption{Computing the region $ \mathcal{S}(1 \%, 1.3)$}
		\label{fig:threshold_regions1} 

\end{center}
\end{figure}
	
\begin{figure}[h!]
\begin{center}	
	
		\subfloat[Proposed]{
			\includegraphics{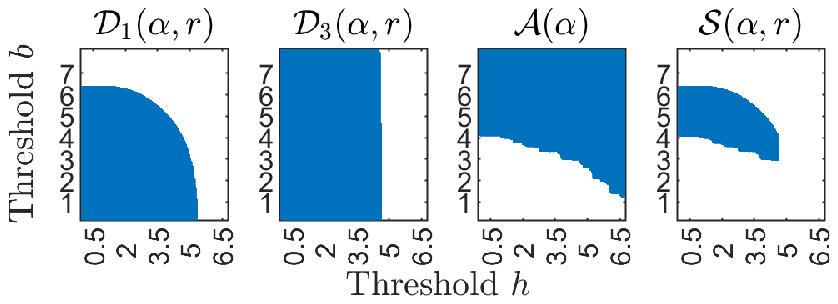}  
			\label{fig:proposed_threshold_regions_2}
		}
		
		\subfloat[Matrix CuSum]{
			\includegraphics{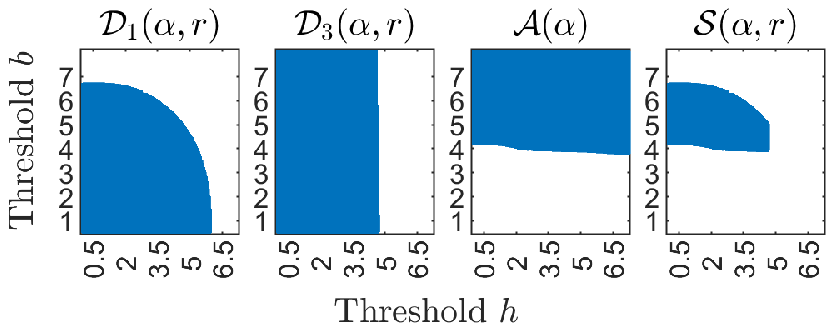}  
			\label{fig:matrix_cusum_threshold_regions_2}
		}

		\captionsetup{justification=raggedright, singlelinecheck=false}
		\caption{Computing the region $ \mathcal{S}(1 \%, 2)$}
		\label{fig:threshold_regions2} 
\end{center}
\end{figure}
\end{appendices}
\newpage

\bibliographystyle{IEEEtranN}
\bibliography{DetIsoARXIV.bbl}    

\end{document}